\documentclass[preprint,12pt]{elsarticle}

\usepackage[titletoc,toc,title]{appendix}
\usepackage{amstext}
\usepackage{enumerate}
\usepackage{epsf}
\usepackage{psfrag}
\usepackage{wrapfig}
\usepackage{lastpage}
\usepackage{float}
\usepackage[english]{babel}
\usepackage{graphicx}
\usepackage[ansinew, latin1]{inputenc}
\usepackage{type1cm}
\usepackage{array}
\usepackage{color}
\usepackage{colortbl}
\usepackage{subfigure}
\usepackage{titlesec}
\usepackage{verbatim}
\usepackage{algorithmicx}
\usepackage{algpseudocode}
\usepackage{algorithm}
\usepackage[normalem]{ulem}
\usepackage{mathtools}
\usepackage{amsmath,color,rotating,psfrag,amsfonts}
\usepackage[dvips]{epsfig}
\usepackage{subfigure}
\usepackage{float}
\usepackage{amsmath}
\usepackage{amssymb}
\usepackage{amsthm}
\usepackage{accents}
\usepackage{epic}
\usepackage[margin=1cm]{caption}
\usepackage[titletoc,toc,title]{appendix}
\usepackage{bbm,bm}
\usepackage{tikz}

\usepackage{booktabs}
\usepackage{siunitx}
\usepackage{multicol}
\usepackage{multirow}

\usepackage{pgfplots}
\usepackage{adjustbox}
\usepackage{titlesec}
\addto\captionsenglish{\renewcommand{\appendixname}{Appendix }}

\newtheorem{remark}{Remark}

\newcommand {\bc} {\mathbf{c}}
\newcommand {\be} {\mathbf{e}}
\newcommand {\beff} {\mathbf{f}}

\newcommand {\bk} {\mathbf{k}}

\newcommand {\bu} {\mathbf{u}}

\newcommand {\bx} {\mathbf{x}}
\newcommand {\bxi} {\boldsymbol{\xi}}
\newcommand {\bPsi} {\boldsymbol{\Psi}}

\newcommand {\bV} {\mathbf{V}}

\newcommand {\bG} {\mathbf{G}}

\newcommand {\bQ} {\mathbf{Q}}
\newcommand {\bR} {\mathbf{R}}
\newcommand {\bT} {\mathbf{T}}

\newcommand\diff{\mathop{}\!\mathrm{d}}

\journal{Journal of Computational Physics}

\begin{document}

\begin{frontmatter}


\title{Level Set Methods for Stochastic Discontinuity Detection in Nonlinear Problems}
\author{Per Pettersson$^{a}$, Alireza Doostan$^{b}$, Jan Nordstr\"{o}m$^{c}$}
\address{$^{a}$NORCE Norwegian Research Centre,
N-5838 Bergen, Norway \\
$^{b}$Aerospace Engineering Sciences, University of Colorado Boulder, CO 80309, USA\\
$^{c}$Department of Mathematics, Link\"{o}ping University, SE-58183 Link\"{o}ping, Sweden}


\begin{abstract}
Stochastic problems governed by nonlinear conservation laws are challenging due to solution discontinuities in stochastic and physical space. In this paper, we present a level set method to track discontinuities in stochastic space by solving a Hamilton-Jacobi equation. By introducing a speed function that vanishes at discontinuities, the iso-zero of the level set problem coincide with the discontinuities of the conservation law. The level set problem is solved on a sequence of successively finer grids in stochastic space. The method is adaptive in the sense that costly evaluations of the conservation law of interest  are only performed in the vicinity of the discontinuities during the refinement stage. In regions of stochastic space where the solution is smooth, a surrogate method replaces expensive evaluations of the conservation law. The proposed method is tested in conjunction with different sets of localized orthogonal basis functions on simplex elements, as well as frames based on piecewise polynomials conforming to the level set function. The performance of the proposed method is compared to existing adaptive multi-element generalized polynomial chaos methods.
\end{abstract}

\end{frontmatter}

\section{Introduction}
Solutions of nonlinear conservation laws often come with uncertainty, and estimated Quantities of Interest (QI) are therefore unreliable. This can be due to unknown material parameters and lack of knowledge on the exact form of the physical laws describing the problems. Uncertainty quantification can be used to characterize these uncertainties in input parameters, geometry, initial and boundary conditions, and to propagate them through the governing equations to obtain statistics of QI.

For problems where the QI depend smoothly on the uncertain input variables, the generalized Polynomial Chaos (gPC) framework offers a range of efficient methods for uncertainty quantification~\cite{Ghanem_Spanos_91,Xiu_Karniadakis_02}. This has been demonstrated extensively for diffusive problems, c.f., applications to fluid flow~\cite{Xiu_Karniadakis_03} and heat conduction~\cite{Xiu_Karniadakis_03_b}. Depending on smoothness, efficient representation is also possible for moderately high-dimensional problems~\cite{Xiu05a,Nobile_etal_08,Ma09a,Doostan09,Doostan11a,Doostan13a,Peng14}.

In stochastic nonlinear wave propagations problems, the solution is typically discontinuous in both physical and stochastic space~\cite{Chen_etal_05,Chanstrami09}. Localized gPC based on adaptive partitioning of the stochastic domain was introduced as Multi-Element generalized Polynomial Chaos (ME-gPC) in~\cite{Wan_Karniadakis_05}. This method is attractive since knowledge of the location of solution discontinuities is not required. One relies instead on a local measure of variance as a criterion for adaptivity. Other methods for stochastic discontinuous solutions that do not rely on explicit calculation of discontinuity locations include Pad\'{e} approximation of discontinuous functions using rational functions~\cite{Chanstrami09}, and multi-resolution analysis schemes based on multi-wavelet expansions that are robust to discontinuities due to hierarchical localization in stochastic space~\cite{LeMaitre_etal_04,Schiavazzi14,Schiavazzi17}. Iterative methods for computation of improved spectral expansions of non-smooth solutions to successively suppress the Gibbs phenomenon were introduced in~\cite{Poette_Lucor_12}.

Efficient stochastic representation (e.g., spectral expansions) of QI require knowledge of the location of discontinuities in stochastic space. Tracking discontinuities in high-dimensional spaces is a challenging problem and many existing methods are subject to restrictions on the geometry of the discontinuities. A hyperspherical sparse approximation framework to detect discontinuities in high-dimensional spaces was recently introduced in~\cite{Zhang_etal_16}, but is restricted to connected star-convex regions. Methods for stochastic discontinuity detection based on polynomial annihilation techniques were introduced in~\cite{Archibald_etal_09}, and based on Bayesian inversion in~\cite{Sargsyan_etal_12}. Both works subsequently used piecewise gPC for stochastic representation of QI on either side of the detected discontinuities. Polynomial annihilation was also used to initialize functional domain decomposition followed by refinement using machine learning with support vector machines to find discontinuities in QI~\cite{Gorodetsky_Marzouk_14}.
A different approach under the name of Multi-Element Minumum Spanning Trees consists in sampling a QI using a minimum spanning tree algorithm that adaptively concentrates samples in stochastic regions with large QI gradients. The stochastic domain is partitioned into nonoverlapping elements of piecewise smooth QI, where the element boundaries are identified by support vector machines~\cite{vanHalder_etal_18}.

For moderately high-dimensional stochastic problems, a viable option to track discontinuities with complex geometries is offered by level set methods. Since the introduction of numerical methods for the solution of level set problems in the seminal work~\cite{Osher_etal_88}, these methods have received considerable attention and applications in image processing, fluid mechanics and materials science~\cite{Sethian_99}, among others. Level set methods are not restricted to star-convex regions and are therefore of interest for complex discontinuous problems of limited dimensionality~\cite{Malladi_etal_94}. The advantages of using level set methods include:
\begin{itemize}
\item
Complex shape recognition: The level set method can handle complex shapes with, e.g., sharp corners and pinch-offs. The complex shapes do not necessarily need to be connected, allowing multiple discontinuities to be captured.
\item
Noise reduction: If the images -- here, the partial differential equation (PDE) solutions -- are  contaminated by noise, level set methods have shown to be an efficient means for shape recovery~\cite{Malladi_Sethian_96}.

\item The level set function provides an implicit parameterization of the discontinuities and can be used for book-keeping when the QI is evaluated at new points in parameter space since its value denotes the signed distance from the discontinuity.
\end{itemize}
In the context of uncertainty quantification, shape recovery was performed on a set of random images and combined with polynomial chaos representation to identify a suitable random parameterization of uncertain geometries in~\cite{Stefanou_etal_09}. Level set methods have also been used for problems with random geometries in~\cite{Lang_etal_13}, where an extended stochastic finite element method with gPC representation and basis enrichment was proposed. A gPC formulation of level set problems for image segmentation through stochastic Galerkin projection was presented in~\cite{Patz_Preusser_14}.

In this paper, we present a new method to track discontinuities in stochastic space by solving a sequence of successively more refined level set problems, constructed such that their iso-zero coincide with the discontinuities of the conservation law we wish to solve. The level set problems are described by a Hamilton-Jacobi equation with a speed function that vanishes at discontinuities of the conservation law. The method is adaptive in the sense that costly evaluations of the conservation law are only performed in the vicinity of the discontinuities during the refinement stage. To the best of our knowledge, this combined method has not previously been considered in the literature. 

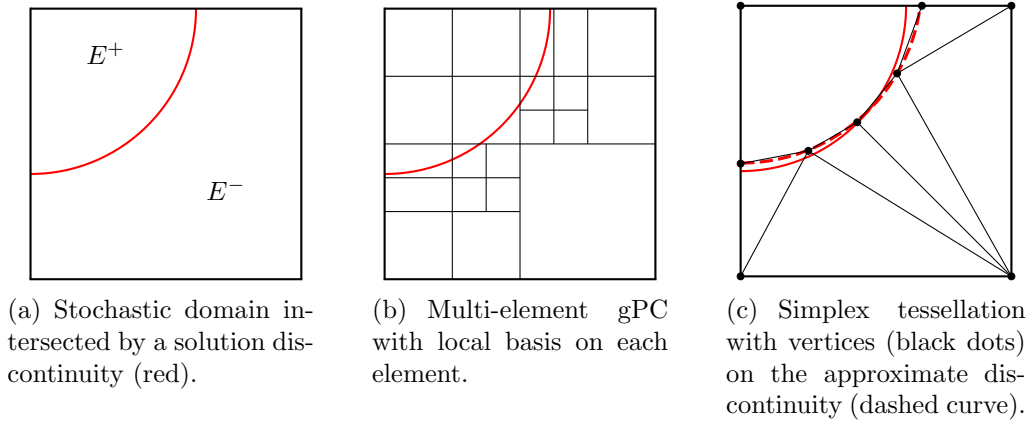
\begin{figure}[H]

  \centering
 \subfigure[Stochastic domain intersected by a solution discontinuity (red).]{
 \label{fig:schematic_a}

 \begin{tikzpicture}[scale=1](10,10)

\draw[thick,red] ([shift=(270:2.2cm)]0,3.6) arc (270:360:2.2cm);

\draw[thick](0,0) -- (0,3.614);
\draw[thick] (-0.014,0) -- (3.614,0);
\draw[thick] (3.6,0) -- (3.6,3.6);
\draw[thick] (0,3.6) -- (3.614,3.6);



\node (E+) at (1, 3) {{\footnotesize$E^{+}$}};

\node (E-) at (2.6, 1.2) {{\footnotesize $E^{-}$}};

\end{tikzpicture}
  }
  \hspace{0.5cm}
 \subfigure[Multi-element gPC with local basis on each element.]{
  \label{fig:schematic_b}
 \begin{tikzpicture}[scale=1](10,10)

\draw[thick,red] ([shift=(270:2.2cm)]0,3.6) arc (270:360:2.2cm);

\draw[thick](0,0) -- (0,3.614);
\draw[thick] (-0.014,0) -- (3.614,0);
\draw[thick] (3.6,0) -- (3.6,3.6);
\draw[thick] (0,3.6) -- (3.614,3.6);

\draw (0.9,0) -- (0.9,3.6);
\draw (1.8,0) -- (1.8,3.6);
\draw (2.7,1.8) -- (2.7,3.6);
\draw (2.25,3.6) -- (2.25,1.8);
\draw (1.35,0.9) -- (1.35,1.8);

\draw (0,2.7) -- (3.6,2.7);
\draw (0,1.8) -- (3.6,1.8);
\draw (0,0.9) -- (1.8,0.9);
\draw (0,1.35) -- (1.8,1.35);
\draw (1.8,2.25) -- (2.7,2.25);

\end{tikzpicture}
 }
%
%
  \hspace{0.5cm}
%
%
 \subfigure[Simplex tessellation with vertices (black dots) on the approximate discontinuity (dashed curve).]{
  \label{fig:schematic_c}
 \begin{tikzpicture}[scale=1](10,10)

\draw[thick,red] ([shift=(270:2.2cm)]0,3.6) arc (270:360:2.2cm);

\draw[line width=1.2pt,red,dash pattern={on 5pt off 2pt} ] ([shift=(265:2.2cm)]0.2,3.7) arc (270:353:2.4cm);

\draw[thick](0,0) -- (0,3.614);
\draw[thick] (-0.014,0) -- (3.614,0);
\draw[thick] (3.6,0) -- (3.6,3.6);
\draw[thick] (0,3.6) -- (3.614,3.6);

\filldraw [black] (0,0) circle (1.3pt);
\filldraw [black] (0,3.6) circle (1.3pt);
\filldraw [black] (3.6,0) circle (1.3pt);
\filldraw [black] (3.6,3.6) circle (1.3pt);

\filldraw [black] (1.55,2.05) circle (1.3pt);
\filldraw [black] (2.41,3.6) circle (1.3pt);
\filldraw [black] (0,1.5) circle (1.3pt);
\filldraw [black] (0.9,1.67) circle (1.3pt);
\filldraw [black] (2.08,2.7) circle (1.3pt);

\draw (0,0) -- (0.9,1.67);
\draw (0,1.5) -- (0.9,1.67);

\draw (3.6,0) -- (0.9,1.67);
\draw (3.6,0) -- (1.55,2.05);
\draw (3.6,0) -- (2.08,2.7);

\draw (3.6,3.6) -- (2.08,2.7);
\draw (2.41,3.6) -- (2.08,2.7);
\draw (1.55,2.05) -- (2.08,2.7);
\draw (1.55,2.05) -- (0.9,1.67);

\end{tikzpicture}
 }


  \caption{Function on stochastic domain divided by solution discontinuity (red curve) and localization of the stochastic surrogate model using frames, ME-gPC, and simpex elements, respectively. The solution is continuous on the subdomains $E^{+}$ and $E^{-}$, respectively, where superscripts $+$ and $-$ refer to the sign of the associated level set function.}
  \label{fig:schematic_of_methods}
\end{figure}

The location of discontinuities in the QI estimated by the level set function are subsequently used to contruct surrogate models from which statistics can be obtained at a low computational cost. Surrogate models based on the gPC framework can be achieved by various means, and we will present several methods in this work. To illustrate the general setting, let $E$ denote the image of a multidimensional random parameterization. Assume that the QI is dependent on a piecewise continuous function on the stochastic subdomains $E^{+}$  and $E^{-}$ as shown in Figure~\ref{fig:schematic_a}. We compare the performance of existing adaptive ME-gPC methods~\cite{Wan_Karniadakis_05} based on partitioning of the stochastic domain into hyper-rectangles (Figure~\ref{fig:schematic_b}). Then we construct localized bases on simplex subdomains obtained from a Delaunay tessellation defined by points on the computed discontinuity and the domain boundaries, as illustrated in Figure~\ref{fig:schematic_c}. We also investigate the performance of \textit{frames} based on piecewise polynomials defined directly by the subdomains $E^{+}$  and $E^{-}$ of Figure~\ref{fig:schematic_a}, where we use the framework in~\cite{Adcock_Huybrechs_16,Adcock_Huybrechs_18}. For all choices of stochastic basis functions, overdetermined systems of equations must be solved to recover the surrogate function and for that purpose we will use $\ell_p$ regression for $p=1,2$.


The paper is organized as follows. The stochastic conservation law is presented in Section~\ref{sec:stoch_cons_laws}. A level set formulation to track discontinuities in the solutions of the stochastic conservation laws is proposed in Section~\ref{sec:imseg_level_sets}. In Section~\ref{sec:param_unc} we review the representation of uncertainty through gPC and its generalization to multi-element gPC by localizing the stochastic basis to elements of hyper-rectangular shape. To handle more complex geometries efficiently, we introduce multi-elements on simplical domains. We next present global orthogonal polynomials restricted to a subdomain of their original support, resulting in a frame instead of an orthogonal basis.  In Section~\ref{sec:num_algor} we present an adaptive algorithm on multiple stochastic grids and a surrogate method to approximate the solution of the conservation law in regions of smoothness.  Section~\ref{sec:lev_set_basis} deals with the computation of spectral surrogate models, including estimation of spectral coefficients using Least-Squares and Least Absolute Deviations methods. Compared to Ordinary Least-Squares methods, Least Absolute Deviations methods are more costly, but also more robust to extreme values, for instance function evaluations on the opposite side of a discontinuity. An algorithm to obtain a simplex tessellation aligned with the zero level set in stochastic space is described in Section~\ref{sec:tessellation_alg}. The proposed methodology is tested in Section~\ref{sec:num_res} and compared to the adaptive multi-element gPC method~\cite{Wan_Karniadakis_05}. Conclusions are drawn in Section~\ref{sec:conc}.


\section{Stochastic nonlinear conservation laws}
\label{sec:stoch_cons_laws}
Let $D \subset \mathbb{R}^{n}$  ($n=1,2,3$) denote the spatial domain with coordinates $\bx$, and $(\Omega, \mathcal{F},\mathbb{P})$ the probability space with sample space $\Omega$, Borel $\sigma$-algebra $\mathcal{F}$, and probability measure $\mathbb{P}$. Consider a random vector parameterization $\bxi=(\xi_1,...,\xi_{d}): \Omega \rightarrow E$ on this probability space, where $E\subset \mathbb{R}^{d}$ and $\xi_i$ ($i=1,\hdots,d$) are independent random variables with bounded range and probability density functions (PDFs), $\rho_{1}(\xi_1),\hdots,\rho_{d}(\xi_d)$. The joint PDF is denoted $\rho(\bxi) = \rho_1(\xi_1)\hdots \rho_{d}(\xi_d)$ and satisfies $\textup{d}\mathbb{P}=\rho(\bxi)\textup{d}\bxi$.

Consider the conservation law on the physical domain $D$ with boundary $\partial D$,
\begin{equation}
\begin{aligned}
\label{eq:cons_lax_ia}
\bu_{t} + \nabla \cdot \beff(\bu) &= \mathbf{0}, \quad \mbox{in } D \times \Omega \times (0, T],\\
L\bu &= \bm{g}, \quad \mbox{in } \partial D \times \Omega \times (0,T]\\
\bu &= \bu_0, \quad \mbox{in } D \times \Omega, \quad t=0,
\end{aligned},
\end{equation}
where $\bu=\bu(t,\bx,\bxi)$ is the solution vector, $\beff$ is the flux function, and $\nabla$ denotes the standard divergence operator in the physical coordinates. $L$ is a boundary operator and $\bu_0(\bx,\bxi)$ the initial solution. The aim of this paper is to efficiently solve~\eqref{eq:cons_lax_ia} by identifying suitable stochastic representations that conform to discontinuities in stochastic space. For simplicity and ease of presentation, we will consider a scalar solution $\bu=u$ for the rest of this paper. 

\begin{remark}
While the proposed strategies could be extended to black box models with input-output relationships, for the interest of discussion, we have chosen to describe the methodology in the context of parameterized, nonlinear PDEs. 
\end{remark}

%
%

\section{Image segmentation for stochastic discontinuity tracking}
\label{sec:imseg_level_sets}

\begin{figure}[H]
\centering
\subfigure[Discontinuous function $u(\bxi)$.]
{\includegraphics[width=0.47\textwidth]{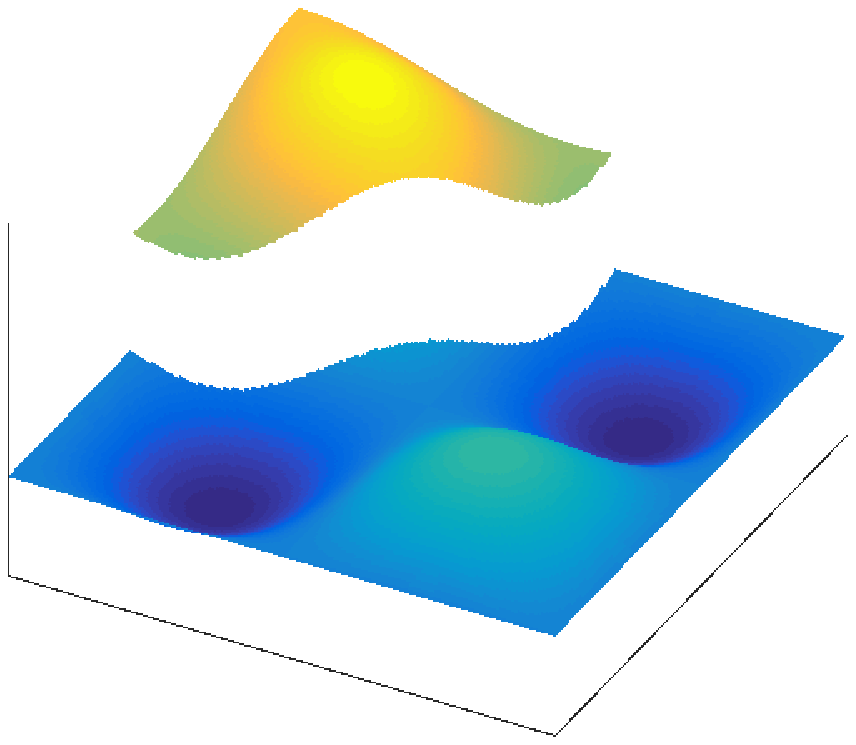}
\label{fig:level_set_schematic_a}
}
\hspace{0.1cm}
\subfigure[Level set function $\phi(\bxi; u)$.]
{\includegraphics[width=0.47\textwidth]{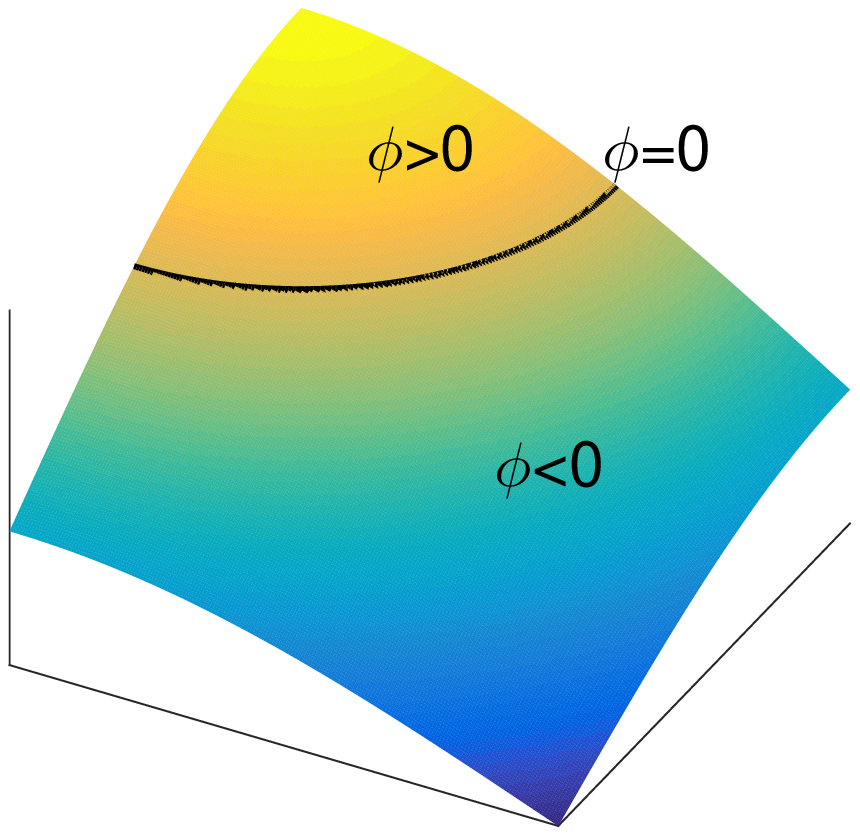}
\label{fig:level_set_schematic_b}
}

\caption{Discontinuous function $u(\bxi)$ (for fixed space and time) and an associated level set function $\phi$ with the zero level set being equal to the location of the discontinuity in $u$.}
	\label{fig:level_set_schematic}
\end{figure}

Our aim is to efficiently solve~\eqref{eq:cons_lax_ia} and to do that we must track the discontinuities of $u$ in $\bxi$. To this end, we introduce the level set function $\phi(\tau,\bx,\bxi)$ where $\bx, \bxi$ have the same meaning as in~\eqref{eq:cons_lax_ia}, and $\tau$ is a pseudo-time. Our goal is that the iso-zero of $\phi$ at some later pseudo-time $\tau$ coincides with the discontinuity location of the solution of \eqref{eq:cons_lax_ia} at some (physical) time of interest $T$. The initial value of the level set function should not necessarily need to coincide with any discontinuity location of $u$. Figure~\ref{fig:level_set_schematic} schematically depicts the properties a level set function should satisfy at large pseudo-times. The discontinuity location of a function $u(\bxi)$ in Figure~\ref{fig:level_set_schematic_a} equal the zero contour of the level set function in Figure~\ref{fig:level_set_schematic_b}. Note that the zero level set function defines a partition of the stochastic domain; solutions to problems where the discontinuities in parameter space do not form closed hypersurfaces or separating hyperplanes, cannot be solved with the proposed level set method. This includes e.g. the Kraichnan-Orszag problem~\cite{Orszag_Bissonnette_67,Wan_Karniadakis_05,Schiavazzi14}, where the solution is continuous but has discontinuous derivatives with respect to the stochastic inputs of the problem.

The evolution of the level set function as introduced in~\cite{Osher_etal_88}, can be described by the PDE
\begin{equation}
\begin{aligned}
\frac{\partial \phi(\tau,\bxi)}{\partial \tau} + F\left(\tau,\bxi \right) \left| \nabla_{\bxi} \phi \left(\tau,\bxi \right) \right| &= 0, \label{eq:level_set} \\
\phi(0,\bxi) &= \phi_{0}(\bxi),
\end{aligned}
\end{equation}
where $F$ is a speed function to be appropriately defined, and $|\cdot|$ denotes the Euclidean distance. The operator $\nabla_{\bxi}$ is the gradient in the stochastic space, i.e., $\nabla_{\bxi} = (\partial_{\xi_1}, \hdots, \partial_{\xi_d})$. In summary,~\eqref{eq:level_set} is a PDE in the stochastic space, to be solved to track solution discontinuities in~\eqref{eq:cons_lax_ia}. Many statistics of interest are restricted to a point in physical space. In the remainder of this paper, we consider statistics at some fixed point in physical space and time. In this case, we solve~\eqref{eq:cons_lax_ia} a number of times in physical space, each time for a different realization of $\bxi$. Then, for fixed space and time, we solve~\eqref{eq:level_set} once. While not the focus of the present study, for spatially and/or temporally dependent QIs, an approach for solving $\phi$ is to solve~\eqref{eq:level_set} independently for each spatial and/or temporal grid. This allows for an embarrassingly parallel solution of $\phi$.

Starting from an initial function $\phi_{0}$, the level set PDE~\eqref{eq:level_set} is evolved in pseudo-time until the iso-zero hits the location of the discontinuity surface in the stochastic space. The key to achieve this result is to make the speed function $F(\tau,\bxi)$ vanish at the discontinuity location. 
The choice for $F$ introduced in~\cite{Malladi_etal_95} is 
\[
F = (1-\epsilon \kappa) \exp\left(- \left| \nabla (G_{\sigma} * u(T,.)) \right|^2 \right),
\]
where $\epsilon > 0$ is a small parameter, 
and the curvature $\kappa = \nabla \cdot (\nabla \phi / \left| \nabla \phi \right| )$ is a common regularization of the level set function. $G_{\sigma}$ is a Gaussian smoothing filter with bandwidth parameter $\sigma$, and the symbol $*$ denotes the convolution operator. The purpose of the Gaussian filter is to remove noise that may otherwise be falsely interpreted as edges or discontinuities~\cite{Malladi_Sethian_96}. Although noise removal is an important feature of level set methods, there is no noise in the solution to the conservation law~\eqref{eq:cons_lax_ia}. Gaussian filtering is therefore not necessary and by choosing $\sigma$ sufficiently small compared to the stochastic grid size, $G_{\sigma}$ becomes the identity operator. In the discrete setting, all derivatives will be finite due to a nonzero grid size. The speed function can be chosen to attain values arbitrarily close to zero in the vicinity of the discontinuities. 

By solving the level set equation~\eqref{eq:level_set} until the zero level set becomes immobile, discontinuities in physical and random space are tracked and this information can be used to construct a local approximation scheme (i.e., adjusted to the discontinuity locations). The steady-state is assumed to be reached when the level set function ceases to change sign at the numerical grid points. This implies that the steady-state zero level set has been found, but non-zero level curves may still change over time.

In this work we assume that we are interested in the solution at some fixed time, but there may be situations where the solution at several different times is desirable. In that case, a discontinuity can emerge at some time instant in a coarsely resolved region where no discontinuity was present at an earlier time instant. In the current implementation, the problem then needs to be solved from scratch. However, with some
modification, the final level set solution at one time instant may provide the initial level set function at another point in time.

The numerical cost of a truly $d$-dimensional problem is $\mathcal{O}(m^{d})$, where $m$ is the number of grid points per dimension in stochastic space. A partial cost reduction to be tried later may be to use so-called narrow-band methods~\cite{Chopp_93} that only include discretization of the regions immediately adjacent to the zero-level set. A more substantial cost reduction may be obtained by identifying lower-dimensional subspaces of the discontinuities, e.g., by using an ANOVA (analysis of variance) decomposition. In a high-dimensional stochastic space, the QI may vary discontinuously in some but not all stochastic dimensions. If $d$ is large, the discontinuity can thus be an object of dimension significantly lower than $d$. The multi-dimensional solution $u(\bxi)$ can be decomposed in terms of effects from each one of the random variables $\xi_{1},\dots, \xi_{d}$ separately, joint effects from each pair of random variables, joint effects from all combinations of triples of random variables, and so on. This leads to the ANOVA decomposition,
\begin{equation}
\label{eq:anova}
u(\bxi) = u^{ \varnothing} + \sum_{s=1}^{d}\sum_{i_1 < \dots < i_s} u^{\{ i_1,\dots,i_s\} }(\xi_{i_1},\dots,\xi_{i_s}),
\end{equation}
where $u^{ \varnothing}$ is the mean value function and the terms $u^{\{ i_1,\dots,i_s\} }$ are orthogonal with respect to the measure of $\bxi$. If higher order interaction terms are small enough to be negligible, or if each of the discontinuities is contained within a subspace of smaller dimension than $d$, the ANOVA decomposition can be truncated, and discontinuities can be tracked in each of the remaining terms separately. Note that the negligible dimensions must be fixed to some constant values (this is known as cut-HDMR (High-Dimensional Model Representations) \cite{Rabitz_etal_99} or anchored ANOVA~\cite{Dick_etal_04}), but the choice of these values is rather arbitrary for the purpose of discontinuity identification, provided that the discontinuities are effectively low-dimensional.

%
%

\section{Spectral expansions in random variables}
\label{sec:param_unc}
\subsection{Generalized polynomial chaos expansion}
Let $\{ \psi_{i_k}(\xi_{k}) \}$ be a univariate orthogonal basis with respect to the weight function $\rho_{k}$, $k=1,...,d$, and let $\bk=(k_1,...,k_{d})$ be a non-negative multi-index. A multivariate global basis of orthogonal functions $\{ \psi_{\bk}(\bxi) : |\bk |< \infty \}$  is constructed through tensorization of univariate basis functions, i.e., the product $\psi_{\bk}=\psi_{k_1}...\psi_{k_{d}}$. Any second-order (finite variance) random function $u(\bxi)$ can then be represented through the generalized Polynomial Chaos (gPC) expansion 
\begin{equation}
\label{eq:gPC_exp}
u(\bxi) \approx u_{ \scriptscriptstyle{\textup{gPC}}}^{\scriptscriptstyle{N}}(\bxi) = \sum_{|\bk|\leq N} c_{\bk}\psi_{\bk}(\bxi), 
\end{equation}
which converges to $u$ in $L_{2,\rho}$ as $N \rightarrow \infty$. The gPC coefficients $c_{\bk}$ are given by the projections of $u(\bxi)$ onto the basis functions, i.e., 
 \[
 c_{\bk} = \frac{\int_{\Omega} u(\bxi) \psi_{\bk}(\bxi) \rho(\bxi) \diff \bxi }{\int_{\Omega} \psi_{\bk}^2(\bxi) \rho(\bxi) \diff \bxi } = \frac{\mathbb{E}(u \psi_{\bk}) }{ \mathbb{E}(\psi^2_{\bk}) },
 \]
 where $\mathbb{E}(\cdot)$ denotes the expectation operator with respect to the PDF $\rho$.


\subsection{Multi-Element generalized polynomial chaos}
\label{sec:me_gpc}
The Multi-Element generalized Polynomial Chaos (ME-gPC) was introduced in~\cite{Wan_Karniadakis_05} and generalized to arbitrary probability measures in~\cite{Wan_Karniadakis_06}. The idea is to partition the random domain into hyper-rectangular elements, and introduce an orthogonal gPC basis with local support on each element. Since the elements are disjoint, basis functions from different elements are orthogonal. Let $\be = (e_1,\dots, e_d) \in \mathbb{N}^{d}$ be a multi-index, where each entry $e_i$ is bounded by some integer $n_i$, and define the element $E_{\be} = E_{e_1} \times \dots \times E_{e_d}$, where $E_{e_i}$ is an open or closed interval within the range of random variable $\xi_i$. The set of multi-elements form a partition of the random space,
\[
E = \bigcup_{\be, e_i \leq n_i} E_{\be}, \quad E_{\be_i}\bigcap E_{\be_j} = \emptyset \mbox{ if } \be_i \neq \be_j.
\]
On each element $E_{\be}$, introduce the local random variable $\bxi_\be$ with the conditional PDF
\[
\rho_{\be}(\bxi_{\be} | \bxi \in E_{\be}) = \prod_{i=1}^{d} \rho_{e_i}(\xi_{e_i} | \xi_i \in E_{e_i}),
\] 
where the univariate conditional PDF on stochastic element $E_{e_i}$ of the i$^{th}$ stochastic dimension is given by
\[
\rho_{e_i}(\xi_{e_i} | \xi_i \in E_{e_i}) =  \frac{\rho_i(\xi_{e_i})}{\mathbb{P}(\xi_i \in E_{e_i})}.
\]
The probability $\mathbb{P}(\xi_i \in E_{e_i})$ is assumed to be positive. Let $\{ \psi_{\be, \bk} \}$ be a set of polynomials on element $\be$ and orthogonal with respect to the conditional PDF. Then the ME-gPC expansion is given by
\[
u(\bxi) \approx u_{ \scriptscriptstyle{\textup{ME-gPC}}}^{\scriptscriptstyle{N}}(\bxi) = \sum_{\be, e_i \leq n_i} \sum_{|\bk|\leq N} c_{\be, \bk}\psi_{\be,\bk}(\bxi),
\]
which is a generalization of \eqref{eq:gPC_exp} to multiple elements.


\subsubsection{Adaptivity Criterion for ME-gPC}
The performance of standard ME-gPC can be improved by adaptivity to regions of sharp variation, e.g., a finer partition of elements in the vicinity of discontinuities. 
In~\cite{Wan_Karniadakis_05}, an adaptive ME-gPC method was developed based on the assumption that if the highest-order gPC coefficients of a multi-element are large in magnitude, then the local variability is not resolved in the current basis. To that end and following~\cite{Wan_Karniadakis_05}, define the element-wise ME-gPC coefficent rate of decay in element $E_{\be}$,
\[
\eta_{\be} = \frac{  \sum_{ |\bk| = N} c_{\be, \bk}^{2} \mathbb{E}(\psi_{\be,\bk}^2) }{ \sum_{0< |\bk|\leq N} c_{\be, \bk}^{2} \mathbb{E}(\psi_{\be,\bk}^2)  },
\] 
which is a measure of the relative contribution of the highest order ME-gPC coefficients to the local variance. The element $E_{\be}$ will be split whenever
\begin{equation}
\eta_{\be}^{\alpha} \mathbb{P}(\bxi \in E_{\be}) \geq \theta_1,
\label{eq:split_crit_1}
\end{equation}
where $0 < \alpha < 1$ and $\theta_1$ are chosen by the user. 

To determine along which dimension to split, the sensitivity of each dimension is evaluated from
\[
r_{\be, i} = \frac{ c_{e_i, N }^{2} \mathbb{E}(\psi_{e_i,N}^2)  }{  \sum_{ |\bk| = N} c_{\be, \bk}^{2} \mathbb{E}(\psi_{\be,\bk}^2)  }.
\]
The random element is split in each dimension that satisfies
\begin{equation}
r_{\be, i} \geq \theta_2 \max_{j=1,...,d} r_{\be, j}, \quad i=1,\dots,d,
\label{eq:split_crit_2}
\end{equation}
where $0 < \theta_2 < 1$ is a design parameter, also chosen by the user.

%
%

\subsection{Multi-element generalized polynomial chaos on simplex shaped elements}
\label{sec:megPC_simplex}

Assuming a finite range of all entries of $\bxi$, the stochastic domain can be partitioned into a set of disjoint simplex elements $\{S_{\be} \}$, analogous to the multi-element partition in Section~\ref{sec:me_gpc}. Analytical expressions for an orthonormal total order $N$ basis,
$
\left\{ \psi_{\bm{\alpha}}  \right\} 
$, 
with the multi-index $\bm{\alpha} \in \mathbb{N}_{0}^{d}, |\bm{\alpha}| \leq N$, 
are given in~\cite{Xu_unknown} for general Dirichlet distributions, i.e., multivariate generalizations of beta distributions. Here we are interested in the uniform probability density function 
 over the simplex since this leads to a more direct relation to more general probability measures through the inverse cumulative distribution function (CDF) method. 

In order to derive orthogonal polynomials on an arbitrary simplex, we first start with orthogonal polynomials on the $d$-dimensional unit simplex $S^{d}$, defined by
\[
S^{d} = \left\{ \bxi \in \mathbb{R}^{d} : \xi_i \geq 0 \mbox{ for } i=1,...,d, \sum_{i=1}^{d}\xi_i \leq 1 \right\}.
\]
Following the notation in~\cite{Xu_unknown}, let 
\[
a_j = 2\sum_{i=j+1}^{d+1}\alpha_i + d - j, \quad j=1,\dots,d,
\]
and
\[
|\bxi_{j-1}| = \left\{
\begin{array}{ll}
\sum_{i=1}^{j-1}\xi_i, & \mbox{if } j>1 \\
0 & \mbox{if } j=1
\end{array}
\right., \qquad
|\bm{\alpha}^{j+1}| = \left\{
\begin{array}{ll}
\sum_{i=j+1}^{d}\alpha_i, & \mbox{if } j<d \\
0 & \mbox{if } j=d
\end{array}
\right..
\]
Let $\left\{ P_{\alpha_j}^{a_j, 0} : j=1,\dots,d,\ |\bm{\alpha}| \leq N \right\}$ be the set of Jacobi polynomials of maximum total order $N$. Then the orthonormal simplex polynomials on $S^d$ are given by
\[
\psi_{\bm{\alpha}} (\bxi) = h_{\bm{\alpha}}^{-1} \prod_{j=1}^{d} (1-|\bxi_{j-1}|)^{\alpha_j} P_{\alpha_j}^{(a_j, 0)}\left(\frac{2\xi_j}{1-|\bxi_{j-1}|} - 1\right),
\]
with the scaling factor
\[
h_{\bm{\alpha}}^{-1}  = \sqrt{ \prod_{j=1}^{d} \frac{ (|\bm{\alpha}^{j}| + |\bm{\alpha}^{j+1}| +d-j) ( 2|\bm{\alpha}^{j}| +d-j +1) }{(d-j+1)(d-j)}}.
\]
The preceeding expressions in this subsection are all special cases of the expressions in Section 5.2 in~\cite{Xu_unknown}.

Next, the orthogonal polynomials on the unit simplex will be generalized to arbitrary simplices. Let the $d+1$ vertices of an arbitrary simplex $s^{d} \in \mathbb{R}^{d}$ be denoted by $\bxi^{i}$, $i=1,...,d+1$. A mapping from a point $\bxi$  in $s^{d}$ to the reference simplex $S^{d}$ with barycentric cordinates $\bm{\lambda} = (\lambda_1,...,\lambda_{d})^{T}$ (excluding the redundant coordinate $\lambda_{d+1}$), is given by
\[
\bm{\lambda} = \bT^{-1} (\bxi - \bxi^{d+1}),
\] 
with the matrix $\bT$ defined by $[\bT]_{i,j} = \bxi_{i}^{j}-\bxi_{i}^{d+1}$. Note that the subscript $i$ denotes random dimension and superscript $j$ denotes vertex index.

The set $\left\{ \psi_{\bm{\alpha}} ( \bm{\lambda}(\bxi)) \right\}$ of orthogonal polynomials is a local basis on the domain restricted to the simplex $s^{d}$. The probability measure of the simplex $s^{d}$ is equal to $|\det(\bT)|/(2^d d!)$. The set of simplex polynomials of all simplices constitute an orthogonal basis on the full stochastic domain $E$.

To distinguish between the localized gPC reconstruction based on simplex elements (to be determined by the solution of level set problems), and standard ME-gPC on hyperrectangular domains, we will use the term Simplex Orthogonal Polynomials (SOP) for the former. Note that there is no difference in cost in the computation of statistics whether it is ME-gPC or SOP.


\subsection{Frames based on restrictions of global orthogonal basis functions}

So far we have considered orthogonal basis functions, on the whole stochastic domain (global basis), hyper-rectangular elements (local basis), and simplex elements (local basis), respectively. We will now consider using global basis functions restricted to a subdomain to be determined by the locations of discontinuities in stochastic space. This construction will lead to a loss of orthogonality and the result is not an \textit{orthogonal basis} but a \textit{frame}, which is a generalization of the concept of basis. Unlike a basis, a frame is redundant, i.e., not linearly independent. Under certain conditions, it still provides good approximation properties of QI. A thorough exposition on frames can be found in~\cite{Christensen_02} and frame approximations on irregular domains are investigated in~\cite{Matthysen_Huybrechs_17,Adcock_Huybrechs_18}. 

Assume that the stochastic domain $E$ (i.e., range of $\bxi$) is partitioned into two disjoint regions $E^{+}$ and $E^{-}$ as shown in Figure~\ref{fig:schematic_a}. The superscripts $+$ and $-$ can be interpreted as being on the 'inside' or 'outside' of a closed-curve discontinuity. The problem setups to be considered in this work always result in a closed curve defined by the discontinuity itself, or by a union of the discontinuity and the boundary of the closed stochastic domain. Let $\{\psi_{\bk}\}$ be a set of global basis functions on $E$, e.g., orthogonal polynomials, and define
\[
\psi_{\bk}^{+} \equiv 
\left\{
\begin{array}{ll}
\psi_{\bk}(\bxi) & \bxi \in E^{+} \\
0 & \bxi \in E^{-} 
\end{array}
\right.,
\quad
\psi_{\bk}^{-} \equiv 
\left\{
\begin{array}{ll}
0 & \bxi \in E^{+} \\
\psi_{\bk}(\bxi) & \bxi \in E^{-} 
\end{array}
\right.
.
\]
The frame $\{ \psi_{\bk}^{+}\} \cup \{ \psi_{\bk}^{-}\}$ is dense in $L_{2,\rho}$. Two functions from the same subdomain ($+$ or $-$) are in general not orthogonal. If any function from the set is removed, the set is no longer a basis for $L_{2,\rho}$. This implies that the frame $\{ \psi_{\bk}^{+}\} \cup \{ \psi_{\bk}^{-}\}$ is a Riesz basis. It satisfies the \textit{frame condition}, i.e., for any $u\in L_{2,\rho}$, $u=\sum_{i \in \mathcal{I} }\tilde{u}_{i}\psi_{i}$ for some index set $\mathcal{I}$ it holds that
\[
A  \left\| u \right\|^{2}_{2} \leq \sum_{i \in \mathcal{I}} |\left\langle u, \psi_i \right\rangle|^{2} \leq B  \left\| u \right\|^{2}_{2}
\]
with the frame bounds $A=\min \lambda(\bG)$ and $B=\max \lambda(\bG)$, where $\bG$ is the Gram matrix, $[\bG]_{i,j}=\left\langle \psi_i, \psi_j \right\rangle$. The weight function of the inner product $\left\langle \cdot ,\cdot \right\rangle$ coincides with the PDF of $\mathbb{E}(\cdot)$.
Since the frame is based on gPC basis functions, we will refer to the method as Frame generalized Polynomial Chaos (F-gPC), with the frame representation
\begin{equation}
\label{eq:gPC_frames}
u_{ \scriptscriptstyle{\textup{F-gPC}}}^{\scriptscriptstyle{N}}(\bxi) = \sum_{|\bk|\leq N} c_{\bk}^{+}\psi_{\bk}^{+}(\bxi) + c_{\bk}^{-}\psi_{\bk}^{-}(\bxi). 
\end{equation}
The idea is that a small number of functions from the Riesz basis will lead to an accurate reconstruction of QI if the support of the basis functions is aligned with solution discontinuities.  The generalization to more than two stochastic subdomains is straightforward. 

Collecting all coefficients $\{ c_{\bk}^{+}\}$ and $\{ c_{\bk}^{-}\}$ in the vector $\tilde{\bc}$, and letting $\mathbf{m}_{\psi}$ be the vector with entries $\mathbb{E}(\psi_{\bk}^{+/-})$ with the same indexing, the mean and variance can be computed from, respectively,
\[
\mathbb{E}(u_{ \scriptscriptstyle{\textup{F-gPC}}}^{\scriptscriptstyle{N}}) = \tilde{\bc}^{\scriptscriptstyle{T}} \mathbf{m}_{\psi}, \quad \mbox{and }\quad \textup{Var}(u_{ \scriptscriptstyle{\textup{F-gPC}}}^{\scriptscriptstyle{N}}) = \tilde{\bc}^{\scriptscriptstyle{T}}\bG\tilde{\bc}-(\tilde{\bc}^{\scriptscriptstyle{T}} \mathbf{m}_{\psi})^2.
\]
The expectations in the entries of $\mathbf{m}_{\psi}$ and $\bG$ have discontinuous integrands, and a robust numerical quadrature rule is needed. The integrals are computed with the trapezoidal rules, using only the already existing function evaluations of the QI that were previously used to track the discontinuity with the level set method.

%
%

\section{Adaptive surrogate level set method for discontinuity tracking}
\label{sec:num_algor}

The computational cost of solving the full problem with moderate stochastic dimensionality is primarily dominated by the extensive cost of solving the conservation law~\eqref{eq:cons_lax_ia} many times for different $\bxi$, as needed to evaluate the speed function $F$. Secondarily, a large contribution to the total cost comes from solving the level set problem~\eqref{eq:level_set}. To address both of these problems, we use an adaptive method and solve a sequence of level set problems with a surrogate method to approximate the speed function. On finer grids, the solution from the coarser grids are used as initial functions. These initial functions are already close to the steady state solution on the fine grids, thus reducing the computational cost compared to solving a fine-grid problem with no previous  estimate of the discontinuity locations. To facilitate the understanding of the proposed method, we first present the elements of the algorithm in some detail. The method is then more succinctly summarized in Algorithm~\ref{alg:adapt_surr_disc_track}, where the numbering corresponds to the numbering in the more detailed description below.

\begin{figure}[H]

  \centering
 {\includegraphics[width=0.96\textwidth]{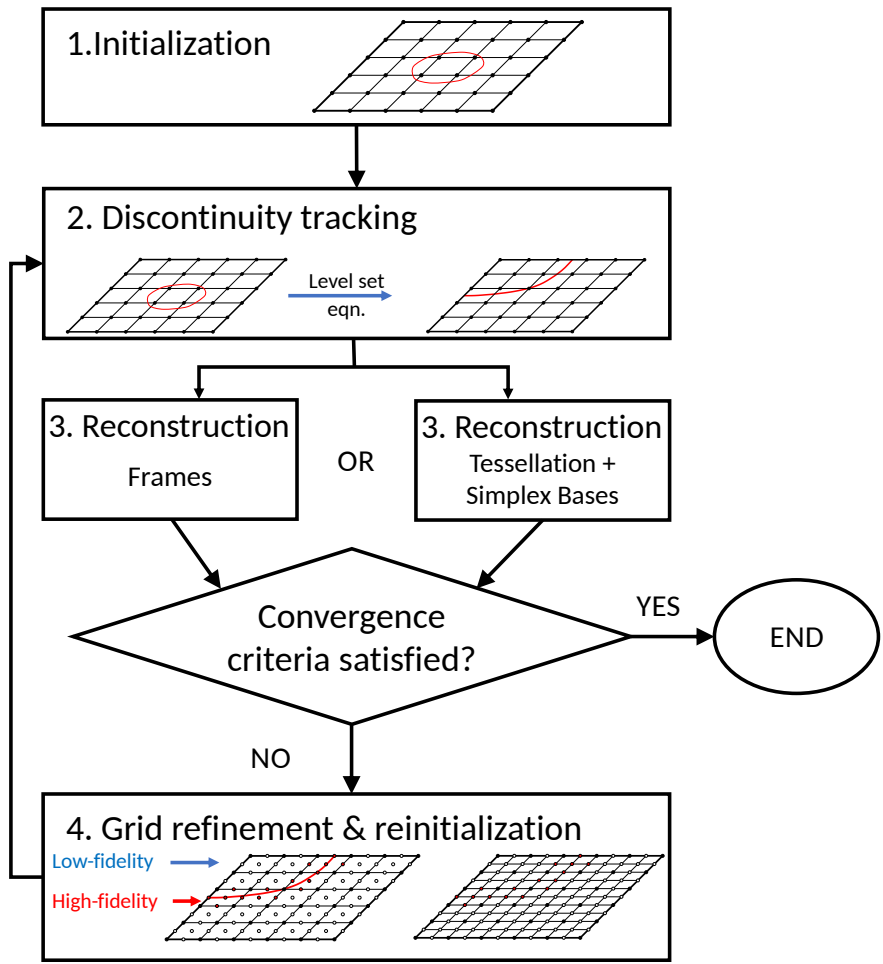}}
 \caption{Graphical representation of Algorithm~\ref{alg:adapt_surr_disc_track}. High-fidelity conservation law solutions are computed in the vicinity of the computed zero level set, and low-fidelity solutions away from the zero level set.}
 \label{fig:flowchart_algorithm}
 
 \end{figure}

\begin{enumerate}
\item
First, the speed function is evaluated on a coarse equidistant grid in $E$ by solving the conservation law~\eqref{eq:cons_lax_ia} once for each stochastic grid point. The level set function is initialized as a small closed curve in the middle of the domain.

\item
The level set problem is solved forward in pseudo-time until locking occurs. The iso-zero of the level set function $\phi$ approximates the location of the discontinuities.

\item
An orthogonal basis or a frame is introduced in stochastic space based on the iso-zero of the level set function. We do one of the following
\begin{itemize}
\item
Simplex tessellation using Delaunay triangulation and computation of local orthogonal basis.
\item
Construction of frames from global orthogonal basis based only on the sign of the level set function (conforming to discontinuities).
\end{itemize}
The solution is reconstructed by solving a Least Absolute Deviations or Ordinary Least Squares regression problem on each element  using the frame/basis functions to be described in Section~\ref{sec:lev_set_basis}. The conservation law evaluations from computing the speed function are used here and no additional solution evaluations are needed. To fit the solution to a local basis of a given simplex element, the number of computed solution points belonging to the element must be larger than the number of local basis functions $(N_{\textup{ev}}>P)$, otherwise the local problem is underdetermined. A simplex element that has no points inside it, most likely has a negligible volume (hence negligible probability) and no basis is introduced on that element. The number of basis functions may vary between elements. 
If the reconstructed solution is sufficiently accurate, as estimated e.g. by checking the decay of the local gPC coefficients, the algorithm terminates. Otherwise, we go to the next step for stochastic grid refinement.

\item
If grid refinement is desired, new evaluations of the conservation law~\eqref{eq:cons_lax_ia} are performed close to the detected discontinuity, as determined by the magnitude of the level set function. In general, this involves numerical solution of~\eqref{eq:cons_lax_ia} on a discrete grid in physical space, even if the QI is defined at some fixed position in space and time.  Away from the discontinuities in $\bxi$-space, the speed function is arbitrary and the solution $u(\bxi)$ is assumed continuous. In these regions, the speed function can be approximated by using a fast proxy method, e.g., evaluating the estimated spectral expansion of the solution on the missing grid points. In the numerical experiments, we simply interpolate new values of $u$ from the ones on the coarse grid. This means that a surrogate method is used to approximate $u$: we use the solution of~\eqref{eq:cons_lax_ia} where high fidelity is needed, and we use interpolation of previous solutions where low fidelity is assumed sufficient. In the test cases of this paper, the interpolation error is negligible compared to the error stemming from estimating the location of discontinuities. 

Next, the level set problem is solved on the refined grid in stochastic space. By using the final solution from the coarse grid as initial function on the refined grid, the number of time-steps until locking occurs can be kept small. The process of grid refinement followed by updated level set solutions can be continued to as many levels of refinement as desired. For clarity, the algorithm is also shown as a graphical flowchart in Figure~\ref{fig:flowchart_algorithm}.

\end{enumerate}


\begin{algorithm}[H]
 \caption{Adaptive Surrogate Method for Discontinuity Tracking.\\ Inputs: $m$ (grid pts per dim.)} 
 {\bf{STEP 1: Initialization of level set problem.}}
\begin{algorithmic} 
\State Assign grid $g^{\textup{coarse}}=\{ \bxi_j\}_{j=1}^{m^d}$.
\State Solve conservation law~\eqref{eq:cons_lax_ia} for $\bxi_{j}$, $j=1,...,m^d$.
\State Initialize $\phi$ and speed function $F$.
\end{algorithmic}
%
{\bf{STEP 2: Track discontinuities.}}
\begin{algorithmic}
\State Solve~\eqref{eq:level_set} in pseudo-time until immobilization of the iso-zero of $\phi$.
\end{algorithmic}
{\bf{STEP 3: Reconstruct stochastic solution.}}
\begin{algorithmic}
\State Tessellation of the domain w.r.t. the iso-zero of $\phi$.
\State Introduce frame/basis functions based on tessellation, assemble overdetermined linear problem.
\State For each element, perform LAD or OLS regression to obtain the element coefficients corresponding to the local frame/basis.
\State Solution meets convergence criteria?
\State \bf{Yes:} Finished. \bf{No:} Go to \bf{STEP 4}.

\end{algorithmic}
{\bf{STEP 4: Grid refinement and reinitialization}}.
\begin{algorithmic} 
\State Refine the grid: add nodes $g^{\textup{new}}$ 
\State $g^{\textup{fine}} \gets g^{\textup{coarse}} \cup g^{\textup{new}} $
\For{$j=1:|g^{\textup{new}}|$}     
\State Evaluate $\phi(\bxi_j)$ by interpolation from surrounding values on $g^{\textup{coarse}}$.
\If{$|\phi(\bxi_j)| < tol$}
\State $u(\bxi_j) \gets$ Conservation law solver.
\Else
\State $u(\bxi_j) \gets$ Surrogate method.
\EndIf
\State Compute $F(u(\bxi_j))$.
\EndFor
\State $g^{\textup{coarse}} \gets g^{\textup{fine}}$
\State Go to {\bf{STEP 2}}.
\end{algorithmic}
\label{alg:adapt_surr_disc_track}
\end{algorithm}

Once the discontinuities in the solution $u$ to~\eqref{eq:cons_lax_ia} have been found by the algorithm presented above, surrogate models for $u$ can be constructed from the spectral expansions introduced in Section~\ref{sec:param_unc}. The topic of the next Section is the computation of the coefficients of the different spectral expansions. In order to avoid additional numerical cost, the surrogate models are exclusively computed from the existing conservation law solutions employed in the level set method.


%
%

\section{ Level set stochastic basis in the presence of discontinuities}
\label{sec:lev_set_basis}
\subsection{Least Absolute Deviations and Ordinary Least Squares Methods}
There are several methods to compute the unknown coefficients $c_{\bk}$ of the ME-gPC, F-gPC, or SOP expansions, e.g., stochastic Galerkin projection~\cite{Ghanem_Spanos_91,Xiu_Karniadakis_02}, non-intrusive spectral projection~\cite{Reagan_etal_03}, stochastic collocation~\cite{Mathelin_etal_05}, and regression~\cite{Berveiller_etal_06}. In this work we use a regression approach, but with non-randomly sampled solutions.
Assuming $N_{\textup{ev}}$ evaluations and $P$ basis functions (i.e., a re-indexing of the multi-index $\bk$ with $|\bk|\leq N$ to single index $k=1,...,P$), we seek a solution to
\[
\left( 
\begin{array}{ccc}
\psi_{1}(\bxi^{\scriptscriptstyle{(1)}}) & \hdots & \psi_{P}(\bxi^{\scriptscriptstyle{(1)}}) \\
\vdots & & \vdots \\
\psi_{1}(\bxi^{\scriptscriptstyle{(N_{\textup{ev}})}}) & \hdots & \psi_{P}(\bxi^{\scriptscriptstyle{(N_{\textup{ev}})}}) 
\end{array}
\right)
\left( 
\begin{array}{c}
c_1 \\
\vdots \\
c_P
\end{array}
\right)
=
\left(
\begin{array}{c}
u(\bxi^{\scriptscriptstyle{(1)}}) \\
\vdots \\
u(\bxi^{\scriptscriptstyle{(N_{\textup{ev}})}})
\end{array}
\right),
\]
i.e., the system $\bPsi \bc = \bu$, where the matrix $\bPsi \in \mathbb{R}^{N_{\textup{ev}}\times P}$ is defined by the entries $[\bPsi]_{i,j} = \psi_{j}(\bxi^{(i)})$, $\bc=\left(c_{1},..., c_{P} \right)^{T}$, and $\bu=\left(u(\bxi^{(1)},..., u(\bxi^{(N_{\textup{ev}})}\right)^{T}$. We choose $P$ so that the system of equations is overdetermined ($P<N_{\textup{ev}}$) and we seek the solution either to the Least Absolute Deviation (LAD) problem
\begin{equation}
\label{eq:LAD}
\min_{\bc} \left\| \bu - \bPsi \bc \right\|_{1},
\end{equation}
or the Ordinary Least Squares (OLS) problem
\begin{equation}
\label{eq:ols}
\min_{\bc} \left\| \bu - \bPsi \bc \right\|_{2}.
\end{equation}
To solve the LAD problem~\eqref{eq:LAD}, we follow~\cite{Shin_Xiu_16,Candes_Tao_05} and perform a pivoted reduced QR factorization of $\bPsi$, i.e., 
\[
\bPsi = \hat{\bQ}\bR, 
\]
where $\hat{\bQ} \in \mathbb{R}^{N_{\textup{ev}} \times r}$ is an orthonormal matrix of rank $r\leq P$ and $\bR \in \mathbb{R}^{r \times P}$ is an upper triangular matrix. Let $\tilde{\bQ} \in \mathbb{R}^{N_{\textup{ev}} \times (N_{\textup{ev}}-r)}$ be an orthonormal matrix with columns orthogonal to those of $\tilde{\bQ}$, i.e., the null space matrix of $\bPsi$. Then we solve the minimization problem
\[
 \min_{\bm{g}} \left\| \bm{g} \right\|_{1} \mbox{ subject to } \tilde{\bQ}^{T}\bm{g} = \tilde{\bQ}^{T}\bu,
\]
with basis pursuit~\cite{vdBerg_Friedlander_09}. The LAD solution is obtained by solving
\[
\bPsi \bc_{\scriptscriptstyle{\textup{LAD}}} = \bu-\bm{g},
\]
for $\bc_{\scriptscriptstyle{\textup{LAD}}} $ in the least-squares sense.

The OLS problem \eqref{eq:ols} has an analytical solution $\bc_{\scriptscriptstyle{\textup{OLS}}}=(\bPsi^{T}\bPsi)^{-1}\bPsi^{T}\bu$ provided $\bPsi$ has full rank, but the LAD solution is more robust to outliers and thus attractive here. 
In more details, in this work, conservation law equations from the opposite side of a discontinuity can be interpreted as outliers in a set of pre-discontinuity solution evaluations. This can happen due to inexact discontinuity identification and is illustrated in Figure~\ref{fig:schematic_c}, where there is an error in the computed discontinuity location (dashed red curve) compared to the exact discontinuity location (solid red curve). When frames are used in the numerical experiments, LAD indeed performs better than OLS in estimating the local frame coefficients, despite some conservation law evaluations being assigned to the wrong solution region. However, the solution response surface obtained from these LAD frame coefficients is still about as erroneous as the response surface obtained using OLS. The explanation is that the error still persists in the partition of the stochastic domain itself. A remedy is implemented by checking the residual $\bPsi \bc_{\scriptscriptstyle{\textup{LAD}}}-\bu$, and re-assigning the points in stochastic space where the absolute value of the residual exceeds the magnitude of the minimum jump in solution values over the discontinuities.

In case of ill-conditioning of the matrix $\bPsi$, which is an issue in particular for the case of F-gPC, the OLS method is adapted as suggested in~\cite{Adcock_Huybrechs_16}, i.e., with a singular value decomposition of the (scaled) approximation of the Gram matrix,
\[
 \bPsi^{T}\bPsi = \bV \boldsymbol{\Sigma} \bV^{*},
 \]
 where the columns of the matrix $\bV$ are the left- and right-singular vectors, and $\boldsymbol{\Sigma}$ is the diagonal matrix of singular values $\sigma_n\geq0$, $n=1,\hdots,P$. Let $\varepsilon > 0$ be a tolerance, and $\boldsymbol{\Sigma}^{\varepsilon}$ be the matrix of truncated singular values with $n^{th}$ diagonal entry $\sigma_{n}>\varepsilon$, and 0 otherwise. The modified OLS coefficients are then given by the truncated SVD approximation
\[
\bc_{\scriptscriptstyle{\textup{OLS}}}^{\varepsilon}= \bV(\boldsymbol{\Sigma}^{\varepsilon} )^{\dagger} \bV^{*}\bPsi^{T}\bu, 
\]
where $\dagger$ denotes the Moore-Penrose pseudoinverse. Note that the estimators $\bc_{\scriptscriptstyle{\textup{OLS}}}$ and $\bc_{\scriptscriptstyle{\textup{OLS}}}^{\varepsilon}$ coincide whenever $\bPsi$ has full rank and $\varepsilon=0$. The same tolerance $\varepsilon=1\cdot 10^{-8}$ is used for both the OLS method and the QR factorization of LAD in the numerical experiments. Approximation of the stochastic solution using LAD and OLS for ME-gPC, SOP, and F-gPC, will be compared in Section~\ref{sec:num_res}.

\begin{remark}
The computation of the Simplex and Frame gPC expansions does not require generating additional PDE samples of~\eqref{eq:cons_lax_ia} as the ones used for the level set construction are readily used.
\end{remark}


\subsection{Simplex tessellation}
\label{sec:tessellation_alg}
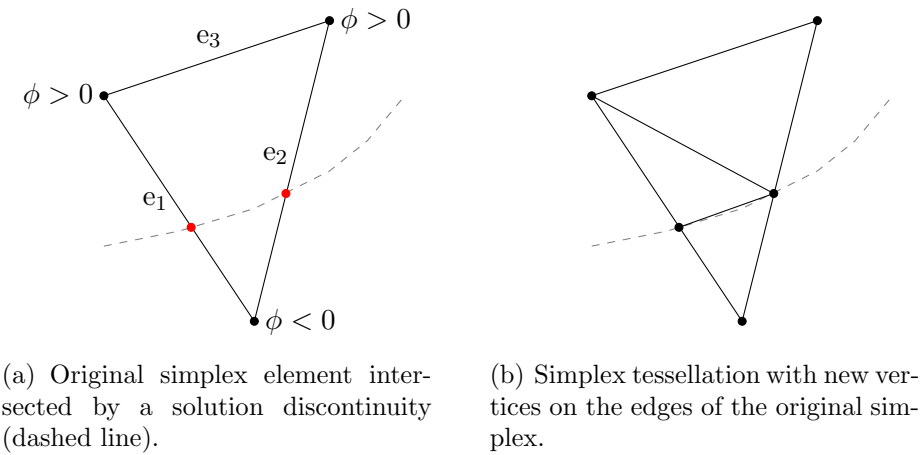
\begin{figure}[H]

  \centering
 \subfigure[Original simplex element intersected by a solution discontinuity (dashed line).]{
 \label{fig:simp_tess_a}

 \begin{tikzpicture}[scale=1](10,10)

\draw (1,4) -- (3,1);
\draw (1,4) -- (4,5);
\draw (4,5) -- (3,1);

\draw[dashed,color=gray] (1,2)--(2,2.2)--(3,2.5)--(4,3.0)--(4.5,3.4)--(5,4);

\filldraw [black] (1,4) circle (1.5pt);
\filldraw [black] (3,1) circle (1.5pt);
\filldraw [black] (4,5) circle (1.5pt);

\filldraw [red] (2.16,2.25) circle (1.5pt);
\filldraw [red] (3.42,2.7) circle (1.5pt);

 \node [left] at (2.0,2.6) {\small{e$_1$}}; 
 \node [left] at (3.6,3.2) {\small{e$_2$}}; 
\node [above] at (2.4,4.5) {\small{e$_3$}}; 

\node [right] at (3,1) {\small{$\phi<0$}};
\node [left] at (1,4) {\small{$\phi>0$}};
\node [right] at (4,5) {\small{$\phi>0$}};

\end{tikzpicture}
  }
  \hspace{0.5cm}
 \subfigure[Simplex tessellation with new vertices on the edges of the original simplex.]{
 \label{fig:simp_tess_b}
 \begin{tikzpicture}[scale=1](10,10)
\draw (1,4) -- (3,1);
\draw (1,4) -- (4,5);
\draw (4,5) -- (3,1);

\draw[dashed,color=gray] (1,2)--(2,2.2)--(3,2.5)--(4,3.0)--(4.5,3.4)--(5,4);

\filldraw [black] (1,4) circle (1.5pt);
\filldraw [black] (3,1) circle (1.5pt);
\filldraw [black] (4,5) circle (1.5pt);

\filldraw [black] (2.16,2.25) circle (1.5pt);
\filldraw [black] (3.42,2.7) circle (1.5pt);

\draw (1,4) -- (3.42,2.7);
\draw (2.16,2.25) -- (3.42,2.7);

\node [right,white] at (3,1) {\small{$\phi<0$}};
\node [left,white] at (1,4) {\small{$\phi>0$}};
\node [right,white] at (4,5) {\small{$\phi>0$}};

\end{tikzpicture}
 }
  \caption{Refinement of simplex tessellation by adding vertices along edges where the level set function $\phi$ changes sign.}
  \label{fig:simp_tess}
\end{figure}

In order to employ orthogonal polynomials restricted to simplex shaped elements, the stochastic domain must be partitioned into a suitable simplex tessellation.
The level set function will be used to iteratively construct a simplex tessellation aligned with the stochastic discontinuities of the QI. In $d$ dimensions, each simplex has $\binom{d+1}{2}$ edges and $d$ vertices. One may directly identify simplex edges that are intersected by the discontinuity by checking if the level set function evaluated at the two vertices defining the edge differ in sign.\footnote{The case of the level set function at the two vertices of an edge being of equal sign implies that the edge is intersected by a discontinuity an even number of times: $0,\ 2,\ \dots$. In this work this fact will be ignored and we assume in this case that the edge is not intersected by any discontinuity. Alternatively, this problem can be remedied by the fact that the computed level set function is known in each grid point. The simplex tessellation can be evaluated qualitatively by checking that all or most level set points within a simplex have the same sign. If this is not the case, for example due to multiple discontinuity intersections, the tessellation should be refined.} The procedure is illustrated for $d=2$ in Figure~\ref{fig:simp_tess}. A simplex element is shown in Figure~\ref{fig:simp_tess_a} with the level set function being positive at two vertices and negative at one. The thin dashed line shows the location of the discontinuity. Due to continuity of the level set function, it must change sign along the edges denoted e$_1$ and e$_2$. The level set function is a distance function and the location along the edges where $\phi=0$ (red dots in Figure~\ref{fig:simp_tess_a}) can be estimated by evaluating a linear relation based on the known values of $\phi$ at the edges. New vertices are added at these points and then a new tessellation is determined from Delaunay triangulation of the updated set of vertices as illustrated in Figure~\ref{fig:simp_tess_b}. The procedure can be iterated, starting from a structured grid of simplices, until convergence of the tessellation. A stopping criterion is enforced by a tolerance on the minimum distance between the vertices is prescribed depending on the spacing between the points of the level set grid. New vertices are introduced only if their distances to existing vertices exceed the tolerance.

%
%

\section{Numerical results}
\label{sec:num_res}

The level set problem~\eqref{eq:level_set} is discretized in stochastic space and pseudo-time by routines from the level set toolbox developed by Ian Mitchell~\cite{Mitchell_08} and modified to the stochastic setting. ENO, WENO and upwind methods are used for the spatial discretization on Cartesian grids, and Runge-Kutta methods for the pseudo-temporal integration. The curvature $\kappa$ is discretized using a second order central discretization, with $\epsilon = 2\Delta \xi$, where $\Delta \xi$ is the grid size of the discretization in the stochastic space. More details on the curvature discretization can be found in the user manual~\cite{Mitchell_07}.

Numerical errors are introduced in the approximation of the discontinuity locations obtained by solving~\eqref{eq:level_set}, as well as in the reconstruction of the solution using either simplex tessellation or estimation of frame coefficients using regression. In order to limit the sources of numerical error and to distinguish between the different errors, we mostly consider numerical test cases with analytical or semi-analytical solutions. This allows for comparison with the exact zero level set function and there is no numerical discretization error in the solution of the conservation law~\eqref{eq:cons_lax_ia}. Results with numerical solutions of the conservation law are included for comparison in case of frame reconstruction. The other methods are more robust and less accurate, so the effect of using numerical conservation law solvers for these cases are expected to be smaller. The range of $\bxi$ is assumed to be $E=[-1,1]^{d}$ with each $\xi_k$ an independent uniform random variable. We use a total-order basis construction which leads to the number of basis functions per element (or number of frame functions for respectively the regions where $\phi>0$ and $\phi<0$) given by $P=(N+d)!/(N!d!)$.

In order to study the performance of the proposed methods, we introduce the discrete $\ell_{1}$ relative error in the reconstruction of the stochastic solution,
\[
\epsilon^{\scriptscriptstyle{\textup{M}}}_{\ell_{1}} = \frac{\left\| u_{\scriptscriptstyle{\textup{M}}} - u_{\textup{ref}} \right\|_{\rho,\ell_1}}{ \left\| u_{\textup{ref}} \right\|_{\rho,\ell_1} }, \quad \mbox{where } \left\| u \right\|_{\rho,\ell_1} \coloneqq  \Delta \xi_1 \dots \Delta \xi_d\sum_{i=1}^{m_{\xi_1}\dots m_{\xi_d}} \rho(\bxi^{(i)}) |u(\bxi^{(i)})|,
\]
for the methods $\textup{M} = \textup{ME-gPC}, \textup{SOP}, \textup{F-gPC}$. The subscript ${\textup{ref}}$ denotes the reference solution, i.e., the exact solution at $\bxi^{(i)}=(\xi_{1}^{(i_1)},\dots,\xi_{d}^{(i_d)})$, and $\Delta \xi_{k}=2/(m_{\xi_k}-1)$ (for $k=1,\dots,d$) is the stochastic grid size. In addition, we present the relative error in means and standard deviations of the solutions,
\[
\epsilon^{\scriptscriptstyle{\textup{M}}}_{\mu} = \left| \frac{ \mu_{\scriptscriptstyle{\textup{M}}} - \mu_{\scriptscriptstyle{\textup{MC}}} }{\mu_{\scriptscriptstyle{\textup{MC}}} } \right|, \quad \epsilon^{\scriptscriptstyle{\textup{M}}}_{\sigma} = \left| \frac{ \sigma_{\scriptscriptstyle{\textup{M}}} - \sigma_{\scriptscriptstyle{\textup{MC}}} }{\sigma_{\scriptscriptstyle{\textup{MC}}} } \right|, 
\]
where $\mu$ and $\sigma$ denote estimators of the mean and standard deviation of the solution $u$, respectively. The subscript $\textup{MC}$ denotes a Monte Carlo reference solution. To reach an error significantly smaller than the errors of the computed solutions using the presented methods, $5\cdot10^{10}$ samples were necessary for the 2D reference solutions, and $5\cdot10^{9}$ samples for the 3D reference solutions.

\subsection{Example 1: Burgers' equation}
Consider the conservation law~\eqref{eq:cons_lax_ia} in two stochastic dimensions and physical domain $D=(-1,1)$, with flux function $f(u)=u^2/2$ and the stochastic Riemann initial condition
\[
u(0,x,\xi_1,\xi_2) = \left\{
\begin{array}{ll}
u_L=a+\sigma_L \cos(c\xi_1) & x \leq x_0,\\
u_R=b+\sigma_R \cos(c \xi_2) & x > x_0,
\end{array}
 \right.
\]
with $a=-b=0.5$, $\sigma_L=0.4$, $\sigma_R=0.3$, $c=3$, $x_0=0$. We use the proposed adaptive level set algorithm on respectively two, three, and four grid levels, always starting on a coarse grid of $m_{\xi_1}=m_{\xi_2}=31$ points. The finest grid will then contain $m_{\xi_1}=m_{\xi_2}=61$ (two levels), $m_{\xi_1}=m_{\xi_2}=121$ (three levels), and $m_{\xi_1}=m_{\xi_2}=241$ (four levels) points in each stochastic coordinate direction.
             
\subsubsection{ME-gPC solution of Burgers' equation}                    
Burgers' equation is solved with the adaptive ME-gPC method for increasing resolution as determined by the two refinement criteria~\eqref{eq:split_crit_1} and~\eqref{eq:split_crit_2} through the choice of the parameters $\theta_1$ and $\theta_2$. The maximum total order of the basis functions of each multi-element is $N=2$, and the results are shown in Table~\ref{tab:conv_adapt_MEgPC} for decreasing tolerance $\theta_1$, and $\theta_2=0.2$.

\begin{table}[h]
\begin{center}
\begin{tabular}{crrccc} \toprule
    {$\theta_1$} & {$|E_{\be}|$} & {$N_{\textup{ev}}$}  & {$\epsilon^{\scriptscriptstyle{\textup{ME-gPC}}}_{\ell_{1}}$}  & {$\epsilon^{\scriptscriptstyle{\textup{ME-gPC}}}_{\mu}$}  & {$\epsilon^{\scriptscriptstyle{\textup{ME-gPC}}}_{\sigma}$}  \\ \midrule
    \midrule
    {$1\cdot 10^{-3}$}  & {144}  & {693}  & 1.20e-1 & 3.55e-3 & 1.05e-2 \\
    {$1\cdot 10^{-4}$}  & {552}  & {2629} & 5.99e-2 & 1.81e-3 & 1.57e-3\\
    {$1\cdot 10^{-5}$} & {1680} & {7933} & 3.52e-2 & 7.16e-5 & 8.34e-4\\
    {$1\cdot 10^{-6}$} & {5564} & {26141} & 1.68e-2 & 1.16e-5 & 2.31e-4 \\
     {$1\cdot 10^{-7}$} & {18204} & {85777} & 9.60e-3 & 6.53e-5 & 4.07e-5
    \\ \bottomrule
\end{tabular}
\end{center}
\caption{Numerical convergence of ME-gPC for different refinement parameter $\theta_1$.}
\label{tab:conv_adapt_MEgPC}
\end{table}

\subsubsection{Level set solution of Burgers' equation: Simplex and Frame gPC}                    
                    
\begin{figure}[H]
\centering
\subfigure[Discontinuities at the zero-level set (red), 65 elements.]
{\includegraphics[width=0.31\textwidth]{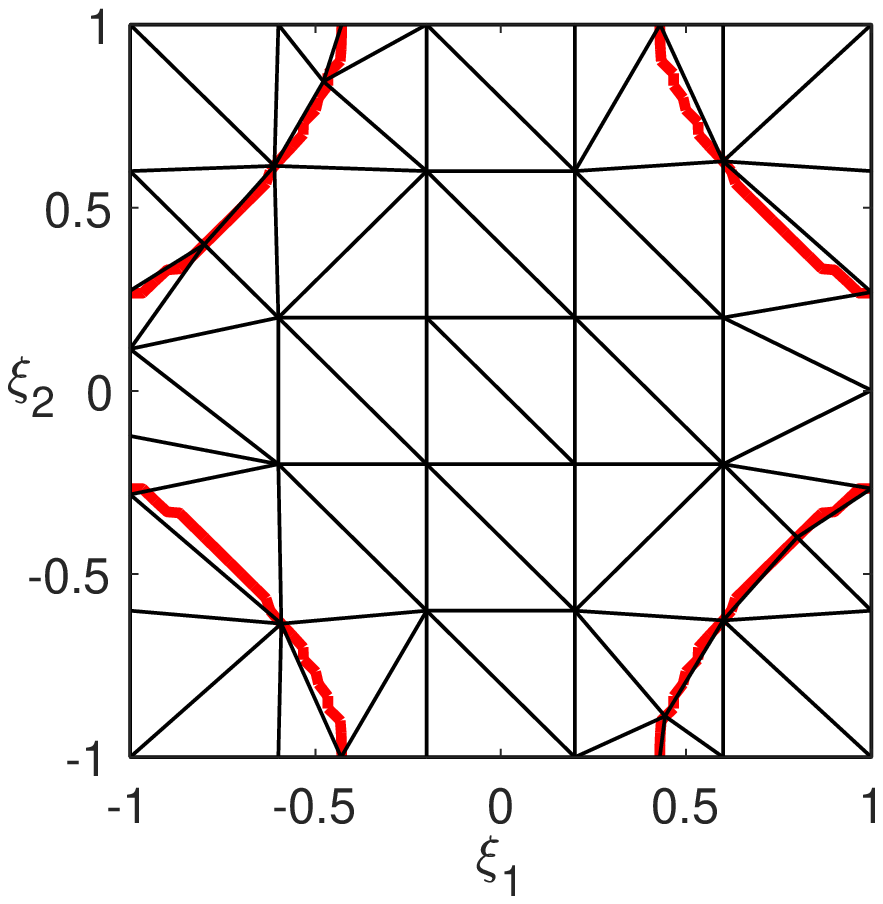}}
\hspace{0.1cm}
\subfigure[Discontinuities at the zero-level set (red), 246 elements.]
{\includegraphics[width=0.31\textwidth]{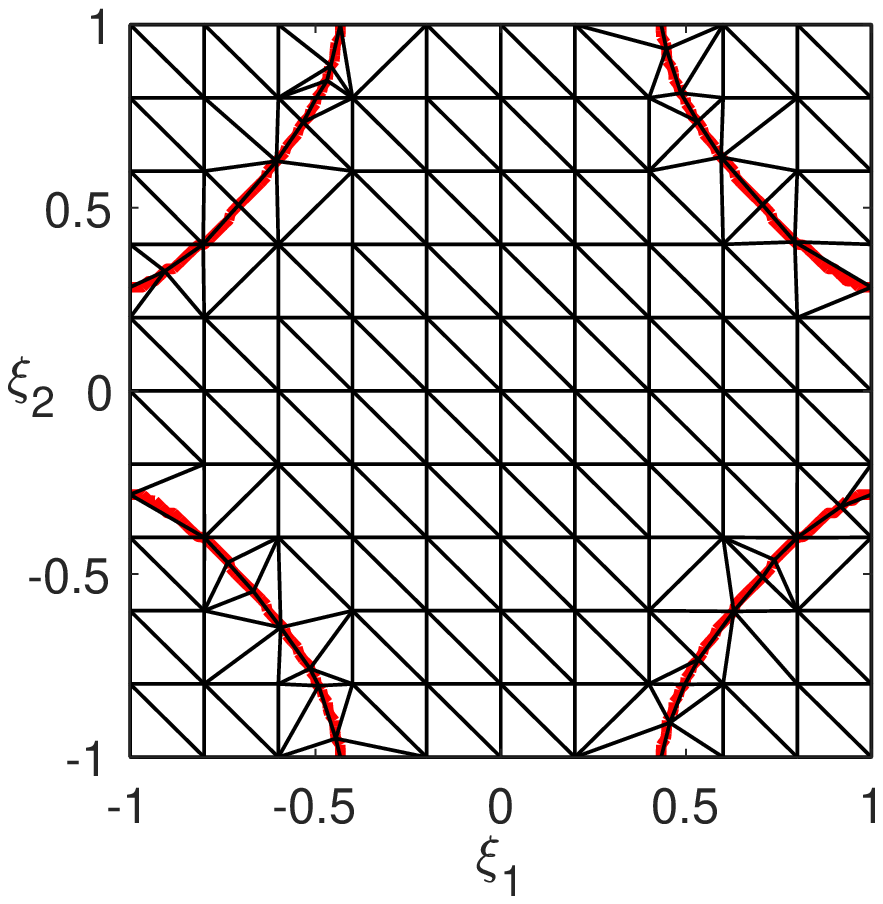}}
\hspace{0.1cm}
\subfigure[Discontinuities at the zero-level set (red), 910 elements.]
{\includegraphics[width=0.31\textwidth]{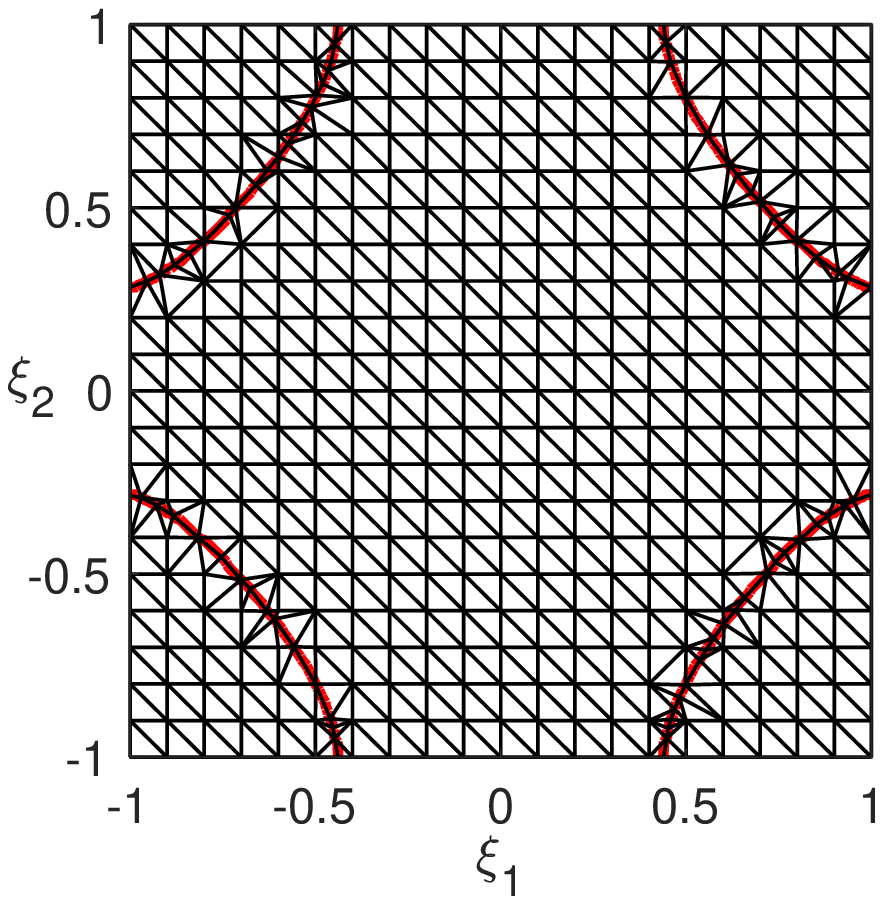}}

\subfigure[High-fidelity (blue) and low-fidelity (red) approx. of $u$.]
{\includegraphics[width=0.31\textwidth]{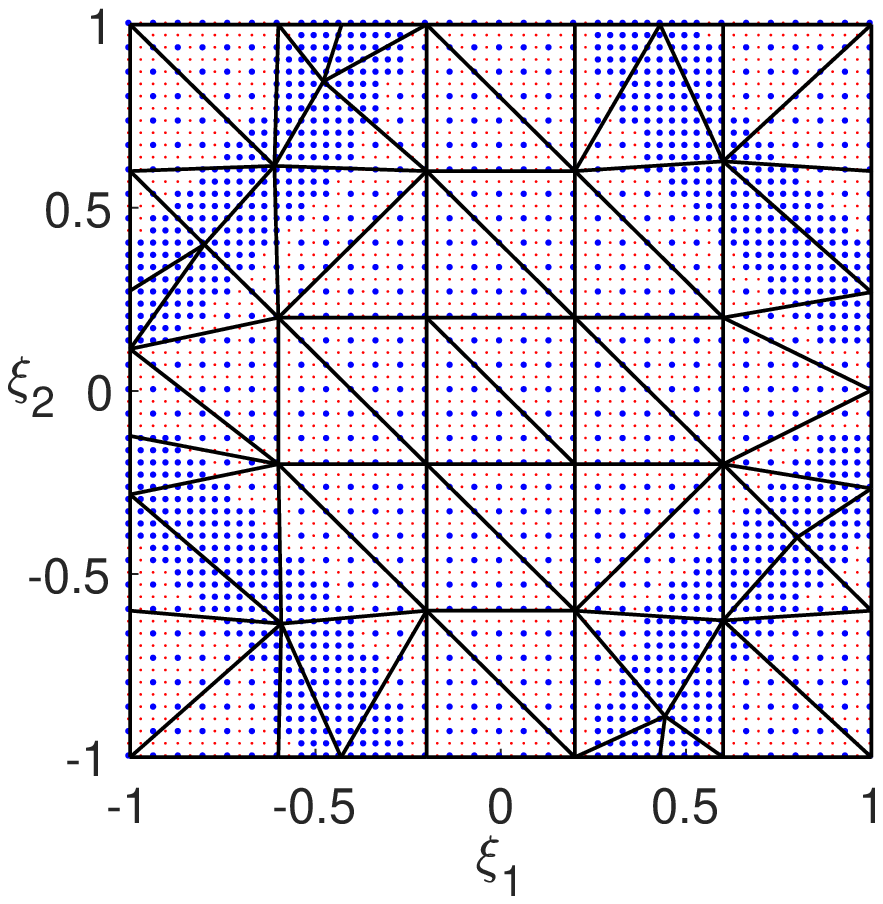}}
\hspace{0.1cm}
\subfigure[High-fidelity (blue) and low-fidelity (red) approx. of $u$.]
{\includegraphics[width=0.31\textwidth]{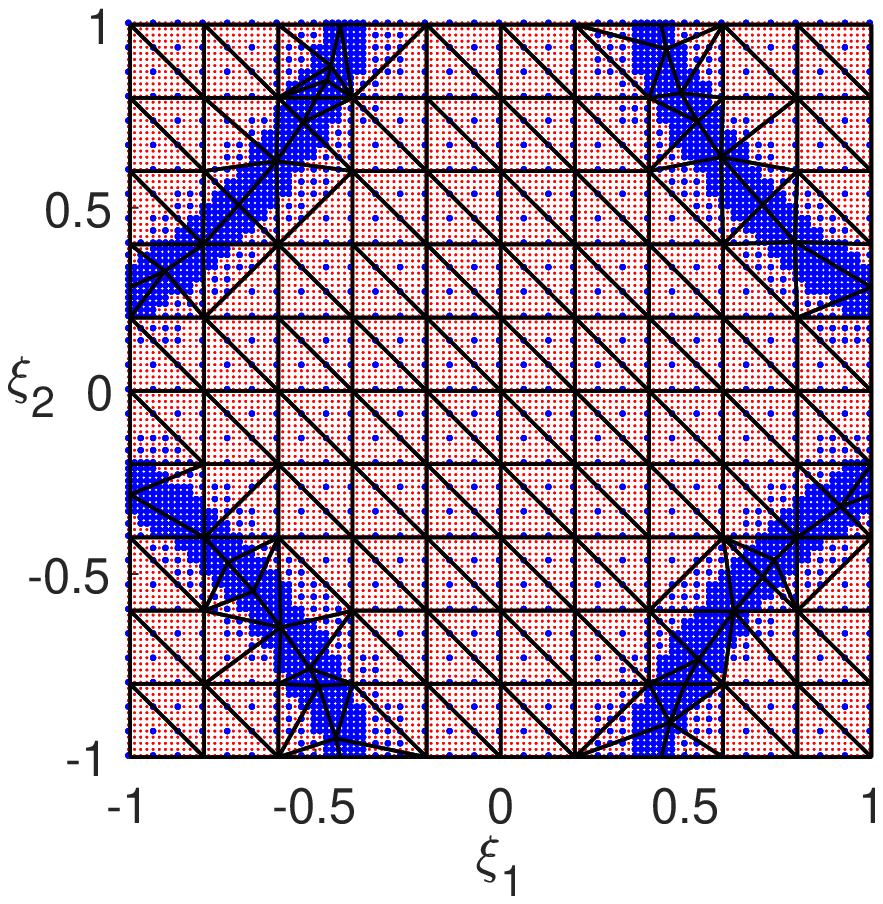}}
\hspace{0.1cm}
\subfigure[High-fidelity (blue) and low-fidelity (red) approx. of $u$.]
{\includegraphics[width=0.31\textwidth]{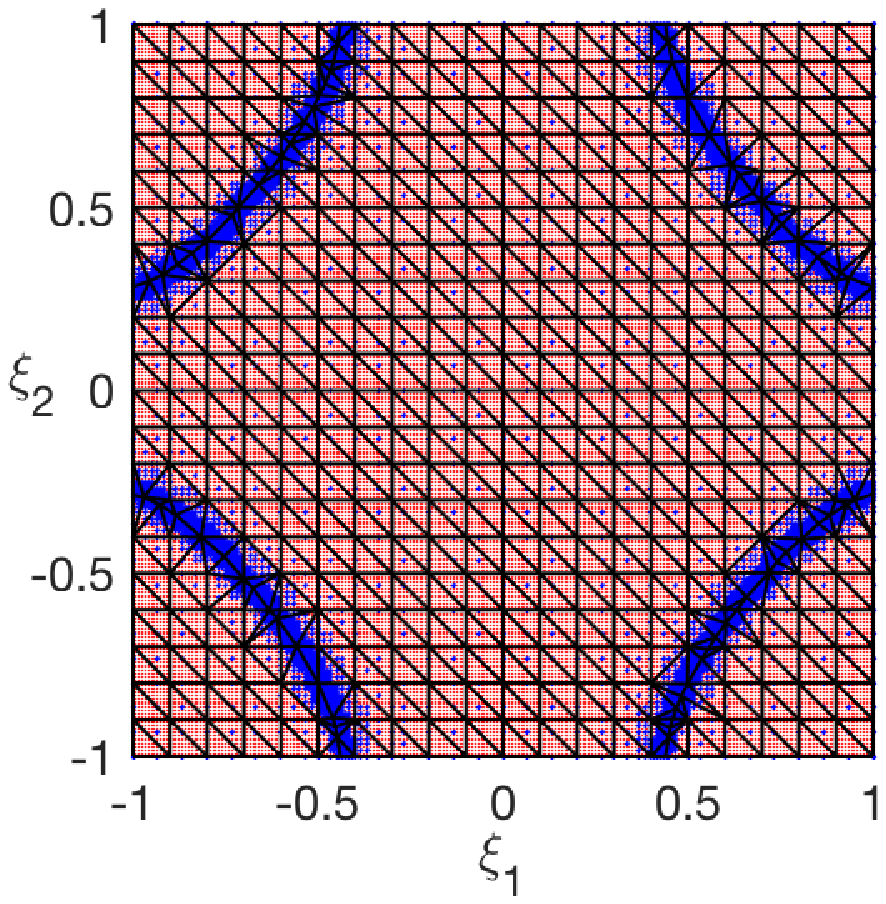}}
\caption{Zero level set contours and conforming simplex tessellations at $x=-0.1$. Two ($m_{\xi_1}=m_{\xi_2}=31,61$), three ($m_{\xi_1}=m_{\xi_2}=31,61,121$) and four grid levels ($m_{\xi_1}=m_{\xi_2}=31,61,121,241$), respectively. The conservation law evaluations are concentrated around the discontinuities.}
	\label{fig:be_234_levels}
\end{figure}

Figure~\ref{fig:be_234_levels} shows the discontinuities identified and the triangulations for the three setups of different number of grids. Note that the triangulations are for reconstruction only; the level set problem has been solved using a sequence of Cartesian grids. The distribution of the high-fidelity conservation law evaluations are indicated by the blue markers in the right figures. They are concentrated around the discontinuities where higher resolution is needed. Red markers indicate grid points where low-fidelity solutions are computed by linear interpolation of the neighboring solutions that were computed on the previous grid level. This leads to computational savings since linear interpolation is significantly cheaper than solving the conservation law. Numerical experiments confirm that the error introduced by replacing the high-fidelity solution with a low-fidelity solution in smooth regions is negligible.

The regression problem \eqref{eq:LAD} is solved using Least Absolute Deviations via basis pursuit within the SPGL1 software~\cite{vdBerg_Friedlander_09}, and the regression problem~\eqref{eq:ols} is solved using OLS. Table~\ref{tab:conv_diff_lev_be_l1} displays the relative $\ell_1$ error for the cases depicted in Figure~\ref{fig:be_234_levels} for different polynomial orders $N$ leading to $P$ basis functions per simplex. As the resolution of the finest grid increases (number of grid points per dimension is denoted $m$), the number of high-fidelity solution evaluations increases, but the proportion $p$ of high-fidelity solutions to the total number of solution estimates (high-fidelity and low-fidelity) decreases, as can be observed in the third column of Table~\ref{tab:conv_diff_lev_be_l1}. On the finer grids, we re-use all high-fidelity function evaluations from the coarser grids, so the number of high-fidelity function evaluations on the finest grid equals the total number of high-fidelity function evaluations. With four grid levels, the numerical cost of the calculation of the speed function is an order of magnitude lower ($p=0.10$) compared to if we had started directly on the finest grid with $m_{\xi_1}=m_{\xi_2}=241$ points per dimension.

\begin{table}[h]
\begin{tabular}{ccc|cc|cc|cc|}
\cline{4-9}
& & & \multicolumn{2}{ c| }{$N=1$, $P=3$} & \multicolumn{2}{ c| }{$N=2$, $P=6$} & \multicolumn{2}{ c| }{$N=3$, $P=10$}  \\ 
\cline{1-9}
\multicolumn{1}{|c}{ $|S_{\be}|$ } & $N_{\textup{ev}}$ & $p$ & LAD & OLS & LAD & OLS & LAD & OLS  \\
\cline{1-9}
\multicolumn{1}{|r}{65} & \multicolumn{1}{r}{1729} & 0.46 & 7.13e-2 & 7.73e-2 & 5.13e-2 & 5.62e-2 & 4.02e-2 & 4.45e-2 \\
\multicolumn{1}{|r}{246} & 3221 & 0.22 & 3.44e-2 & 4.13e-2 & 2.75e-2 & 3.03e-2 & 2.18e-2 & 2.37e-2 \\
\multicolumn{1}{|r}{910} & 5808 & 0.10 & 2.05e-2 & 2.43e-2 & 1.77e-2 & 1.50e-2 & 1.55e-2 & 1.77e-2 \\
\cline{1-9}
\end{tabular}
\caption{Numerical convergence of $\epsilon^{\scriptscriptstyle{\textup{SOP}}}_{\ell_{1}}$ for the Burgers' equation Riemann problem for different orders $N$ of polynomial reconstruction, varying number of simplex elements ($|S_{\be}|$), number of evaluations of the conservation law ($N_{\textup{ev}}$).  The Ordinary Least-Squares (OLS) and Least Absolute Deviations (LAD) methods are used to locally estimate the SOP coefficients. }
\label{tab:conv_diff_lev_be_l1}
\end{table}
The computational cost is furthermore reduced since we start very close to the steady solution at the finest grid levels, which is only possible with coarser grid solutions as initial functions. This is not quantified in the table.

For the finest discretizations where the errors are similar, the ME-gPC method requires more than 4 times as many evaluations of the conservation law solver. Since this is expensive, in particular for computational fluid dynamics problems, the extra cost of solving the level set problem may be negligible.
\begin{remark}
In this work the cost of evaluating the conservation law solution is virtually negligible since we employ analytical solutions to study the performance of the proposed methods and isolate method specific errors without introducing an additional physical discretization error.
\end{remark}

\begin{figure}[H]
\centering
\subfigure[Reconstruction using LAD.]
{\includegraphics[width=0.48\textwidth]{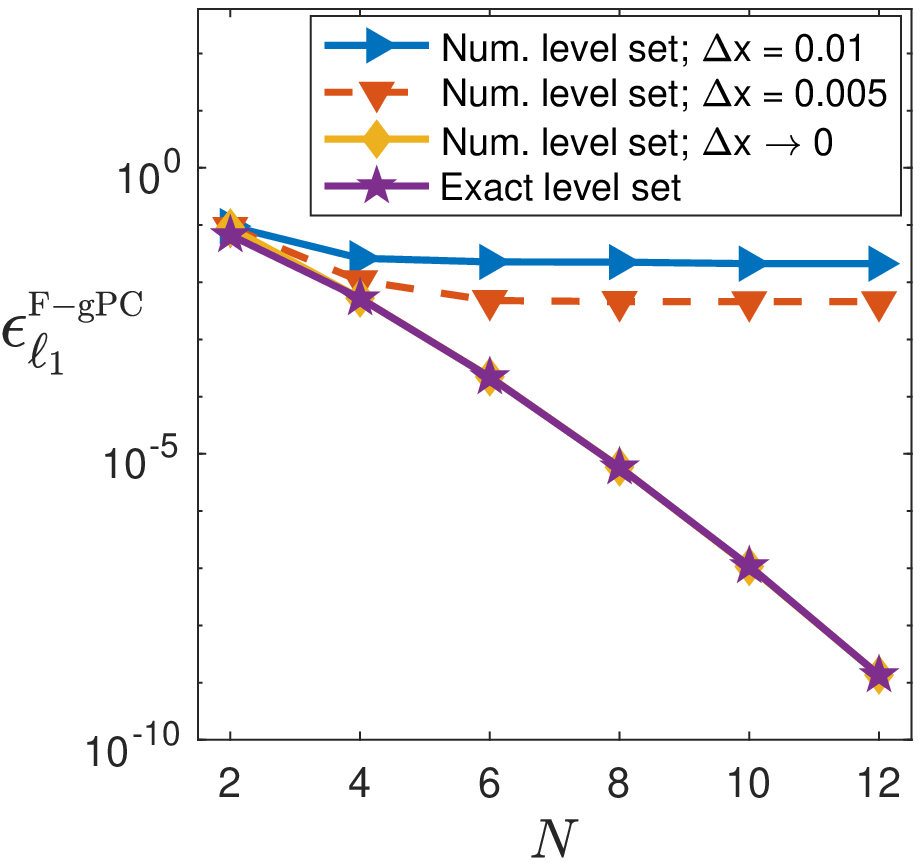}}
\hspace{0.1cm}
\subfigure[Reconstruction using OLS.]
{\includegraphics[width=0.48\textwidth]{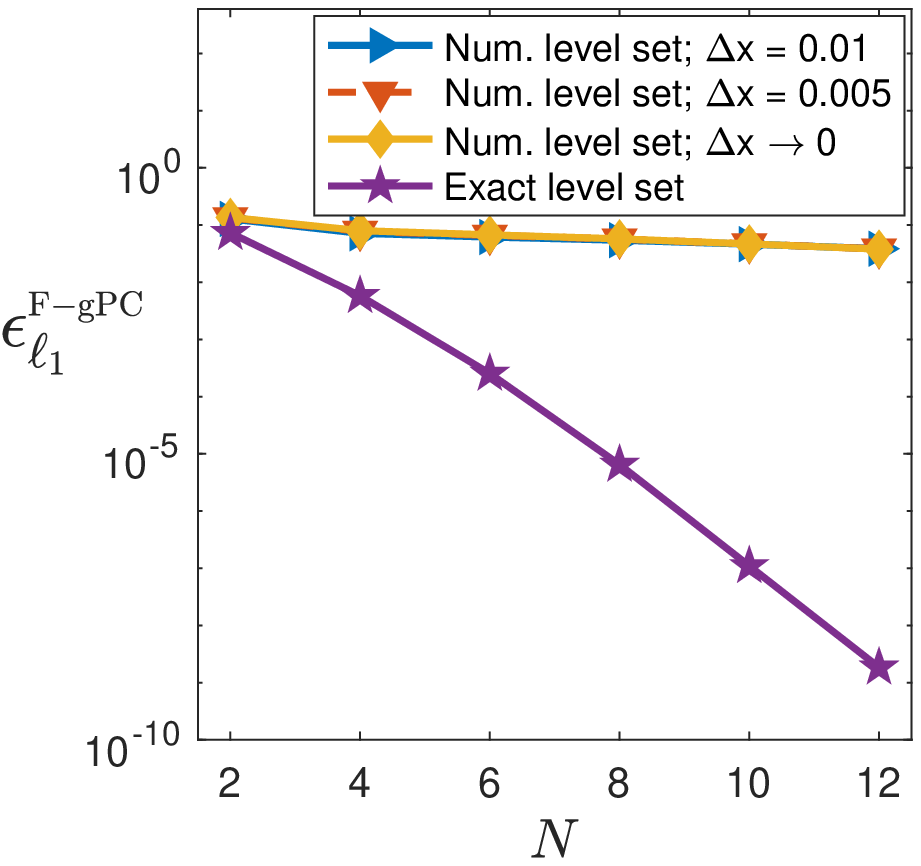}}

\caption{Numerical convergence of $\epsilon^{\scriptscriptstyle{\textup{F-gPC}}}_{\ell_{1}}$ for Burgers' equation for different orders of piecewise polynomial frames, using OLS and LAD. Global total order Legendre polynomials are used for the frames, on a stochastic grid with 61 points per dimension, and $N_{\textup{ev}}=1729$.}
	\label{fig:conv_global_pol_BE_l1_m61}
\end{figure}

Next we use the level set solution for Burgers' equation to define frames, that are piecewise continuous for all $\bxi$ where the level set solution $\phi(\bxi)$ is positive, and negative, respectively. Figure~\ref{fig:conv_global_pol_BE_l1_m61} shows the relative error $\epsilon^{\scriptscriptstyle{\textup{F-gPC}}}_{\ell_{1}}$ for the Burgers' equation test case where the computed frame coefficients have been used to estimate the solution. The LAD and OLS solutions based on the exact zero level set serve as a reference for the error with increasing order of frames. In addition to results based on exact solutions to Burgers' equation, we also present results for a flux-limited Roe scheme~\cite{Roe_81}. Numerical error in the solution to Burgers' equation results in larger relative errors for all polynomial orders $N$. However, for LAD the error decays exponentially if the numerical discretization is sufficiently fine, as evident from the case $\Delta x \rightarrow 0$. In this case we use the analytical solution to Burgers' equation and there is consequently no error in $u$, but we still use the numerical level set solver. As expected, the OLS solution is sensitive to misclassified conservation law solutions, and the relative error never falls behind the order of $10^{-2}$, independently of the accuracy of the numerical solution to Burgers' equation.


\begin{figure}[H]
\centering
\subfigure[Reconstruction using LAD.]
{\includegraphics[width=0.48\textwidth]{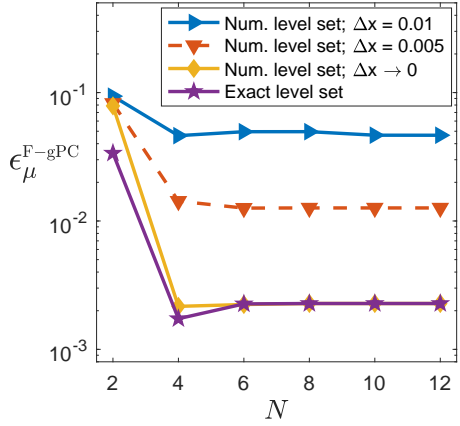}}
\hspace{0.1cm}
\subfigure[Reconstruction using OLS.]
{\includegraphics[width=0.48\textwidth]{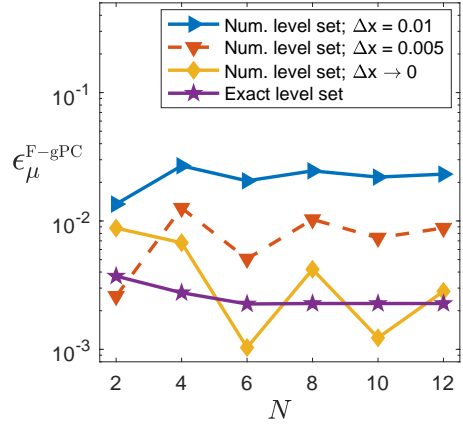}}
\subfigure[Reconstruction using LAD.]
{\includegraphics[width=0.48\textwidth]{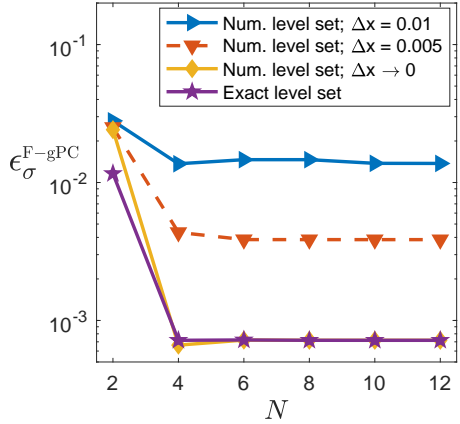}}
\hspace{0.1cm}
\subfigure[Reconstruction using OLS.]
{\includegraphics[width=0.48\textwidth]{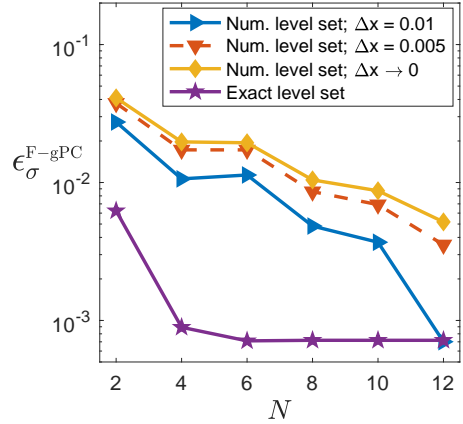}}

\caption{Numerical convergence of $\epsilon^{\scriptscriptstyle{\textup{F-gPC}}}_{\mu}$ and $\epsilon^{\scriptscriptstyle{\textup{F-gPC}}}_{\sigma}$ for Burgers' equation for different orders of polynomial reconstruction, using OLS and LAD. Global total order Legendre polynomials are used for the reconstruction, stochastic grid with 61 points per dimension, and $N_{\textup{ev}}=1729$.}
	\label{fig:conv_global_pol_BE_mu_sigma_m61}
\end{figure}

The error in mean and standard deviation is shown in Figure~\ref{fig:conv_global_pol_BE_mu_sigma_m61}. The trend is less clear than in the case of the error in the $\ell_1$ norm, mainly due to approximation error in the computation of mean and standard deviation using frames.

\subsection{Example 2: CO$_2$ migration model with three stochastic parameters}
\begin{figure}[H]
\begin{center}

\begin{tikzpicture}[scale=1.0]

\draw[thick] (0,0) -- (6,1);
\draw[thick] (0,2) -- (6,3);
\draw[thick] (0,-0.2) -- (6,0.8);

\draw [-latex, black, thick, shorten >= 0.75cm] (6,0.8) -- (7.2,1.0);
\node [right] at (6.3,0.7) {$x$};

\draw[thick] (0,-0.2) -- (6,-0.2);
\node [right] at (1.8, -0.03) {\footnotesize{$\theta$}};
\draw[thin] (1.22,0.02) -- (1.24,-0.2); 
\draw[thin] (1.15,0.008) -- (1.17,-0.2); 


\draw [-latex, black, thin, >= ] (5.70, 1.70) -- (5.495,2.93);
\draw [-latex, black, thin, >= ] (5.70, 1.70) -- (5.835,0.95);
\node [right] at (5.667, 1.95) {\small{$1$}}; 

\draw [-latex, black, thin, >= ] (2.61,2.00) -- (2.76,1.10);
\draw [-latex, black, thin, >= ] (2.61,2.00) -- (2.54,2.42);
\node [left] at (3.11, 1.85) {\small{$u$}};

\node [left] at (2.45, 1.7) {\small{CO$_2$}};

\draw[-latex, thick, >=] (-0.3,0.9) -- (0.3,1.0);
\node[above] at (-0.1,0.95) {\small{$Q$}};

\draw[thick] plot [smooth,tension=0.3] coordinates{(0.3,2.05) (0.6, 1.90) (0.9,1.65) (1.1,1.40) (1.3,1.05) (1.6,0.25)};

\draw[thick] plot [smooth,tension=0.5] coordinates{(2.5,0.4) (2.6,0.75) (2.7,1.0) (3.0,1.4) (3.4,1.8) (4.2,2.35) (5.15,2.85)};
\node [above] at (4.2, 1.35) {\small{brine}}; 

\end{tikzpicture}

\end{center}

\caption{1D vertical equilibrium model of subsurface CO$_2$ transport. The solution $u$ is the relative height of the CO$_2$ plume that migrates due to buoyancy and background flow with rate $Q$.} 
\label{fig:setup_CO2}
\end{figure}
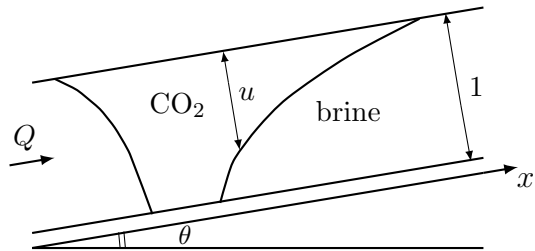

Consider the migration of CO$_2$ in a sloping aquifer in vertical equilibrium with homogeneous material properties, modeled by the hyperbolic nonlinear PDE
\begin{align}
\phi \mathcal{R}\frac{\partial u}{\partial t }+\frac{\partial f(u)}{\partial x} &= 0, \quad \mbox{in } D \times (0, T]. \label{eq:CO2_PDE} \\
u &= u_{0}, \quad t=0,
\end{align}
where $u$ now denotes the normalized height of the CO$_2$ plume, $\phi$ is porosity, and $\mathcal{R}$ is the accumulation coefficient that accounts for trapping of CO$_2$. Additionally,
\[
\mathcal{R} = \left\{ \begin{array}{ll}
1- S_{\textup{br}}-S_{\textup{cr}} & \mbox{if } \frac{\partial u}{\partial t} < 0\\
\\
1- S_{\textup{br}} & \mbox{if } \frac{\partial u}{\partial t} > 0
\end{array}
\right. ,
\]
where $S_{\textup{br}}$ and $S_{\textup{cr}} $ are the residual saturation of brine in CO$_2$, and the residual saturation of CO$_2$ in brine, respectively.
The flux function $f$ is given by
\begin{equation}
\label{eq:frac_flux}
f(u) = \frac{\left( Q+K(1-u) \right)M u}{1+(M-1)u}, 
\end{equation}
with background flow rate $Q$  [L$^{2}$T$^{-1}$], mobility ratio $M$, and $K$ is the product of permeability, density difference between the two phases, buoyancy force, and mobility of brine. The initial time $t=0$ corresponds to the end of the injection period of duration $\tau$, during which CO$_2$ has been injected at a rate $Q_{\text{inject}}$. The initial function $u_0$ describes the shape of the plume during the end of the injection period. More details on the derivation of the model can be found in~\cite{MacMinn_etal_10,Pettersson_16}, including the semi-analytical solution of~\eqref{eq:CO2_PDE}. A sketch of the problem setup is provided in Figure~\ref{fig:setup_CO2}. The slope angle $\theta$ determines the advection speed governed by buoyancy that enters the flux through the parameter $K$. Trapping of CO$_2$ occurs as the plume recedes from a region, driven by buoyancy and background flow, and creates immobile pockets of CO$_2$ left behind in the brine phase.


In the numerical experiments, $D=[0,2000]$, $S_{\textup{br}}=S_{\textup{cr}}=0.1$, $\theta=0.15$, and the initial function $u_0$ is given by Eqn.~(6) in~\cite{Pettersson_16} with $\tau=20$ years. Eqn.~\eqref{eq:CO2_PDE} is solved until $T=600$ (years). We assume that the parameters $M$, $K$, and $Q$ are uncertain. The mobility ratio is given by $M=\lambda_{\textup{c}}/\lambda_{\textup{b}}$, where the CO$_2$ mobility $\lambda_{\textup{c}}$ is uniformly distributed in $ [0.7 \ 1.3] \times 6.25\cdot 10^{-5} $ and $\lambda_{\textup{b}}$ is uniformly distributed in $[0.8, \ 1.2] \times 5\cdot 10^{-4} $. This model takes into account the uncertainty in the endpoints of the relative permeability curves.
The weighted permeability $K$ is lognormal due to the permeability $k$ being lognormal with mean 200 mD and standard deviation 50 mD. The background flow $Q$ is exponentially distributed with mean $1\cdot 10^{-9}$.
The uncertain parameters are represented via the CDF transformations $M=\tilde{F}_{M}^{-1}((\xi_1+1)/2)$, $K=\tilde{F}_{K}^{-1}((\xi_1+1)/2)$, and $Q=\tilde{F}_{Q}^{-1}((\xi_1+1)/2)$, respectively, where $\tilde{F}$ denotes the CDF of each parameter.

\subsubsection{ME-gPC solution of CO$_2$ storage problem}
We solve the CO$_2$ storage problem with piecewise quadratic basis functions on an adaptive ME-gPC discretization of the stochastic space and the results are shown in Table~\ref{tab:conv_adapt_MEgPC_CO2}. The observed convergence rate of the $\ell_1$ error of almost $0.5$ with the number of function evaluations is not surprising. The error in mean and standard deviation shows a similar trend, although the convergence of the mean is not monotone. Compared to standard Monte Carlo simulation, the numerical cost is reduced 2 orders of magnitude. The number of calls to the conservation law solver must be increased significantly to further reduce the error, determined by the splitting parameter $\theta_1$. \begin{table}[htb]
\begin{center}
\begin{tabular}{crrccc} \toprule
    {$\theta_1$} & {$ | E_{\be} |$} & {$N_{\textup{ev}}$}  & {$\epsilon^{\scriptscriptstyle{\textup{ME-gPC}}}_{\ell_{1}}$} & {$\epsilon^{\scriptscriptstyle{\textup{ME-gPC}}}_{\mu}$} & {$\epsilon^{\scriptscriptstyle{\textup{ME-gPC}}}_{\sigma}$}  \\ \midrule
    \midrule
    {$1\cdot 10^{-3}$}  & {346}  & {4285}  & 1.81e-1 & 3.17e-3 & 2.68e-2 \\
    {$1\cdot 10^{-4}$}  & {1603}  & {19211} & 1.08e-1 & 5.04e-3 & 1.33e-2\\
    {$1\cdot 10^{-5}$} & {7761} & {90946} & 4.98e-2 & 7.70e-4 & 6.93e-3\\
    {$1\cdot 10^{-6}$} & {33622} & {386934} & 2.56e-2 & 1.85e-4 & 3.43e-3\\
        {$1\cdot 10^{-7}$} & {161057} & {1832600} & 1.45e-2 & 1.10e-4 & 1.61e-3
    \\ \bottomrule
\end{tabular}
\end{center}
\caption{Numerical convergence of ME-gPC for the CO$_2$ storage transport problem with varying refinement parameter $\theta_1$.}
\label{tab:conv_adapt_MEgPC_CO2}
\end{table}


\subsubsection{Level set solution of CO$_2$ storage problem: Simplex and Frame gPC}
\begin{figure}[H]
\centering
\subfigure[$\xi_3=-1$.]
{\includegraphics[width=0.48\textwidth]{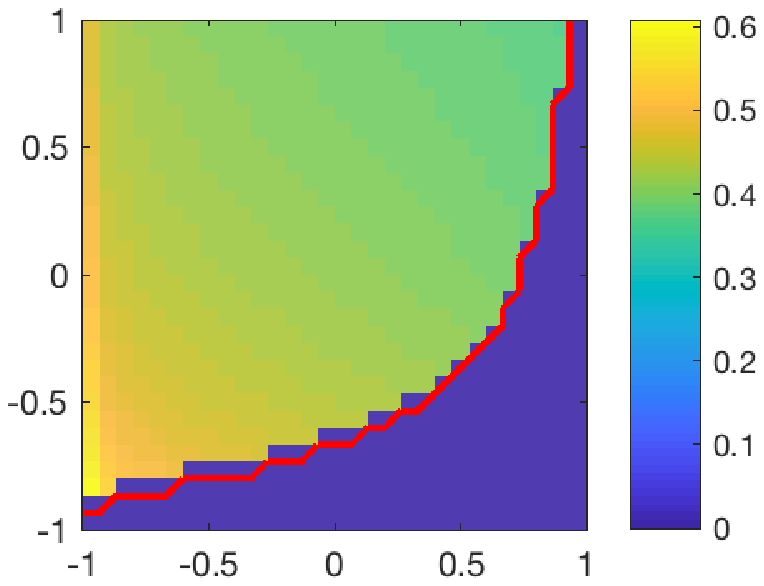}}
\subfigure[$\xi_3=0$.]
{\includegraphics[width=0.48\textwidth]{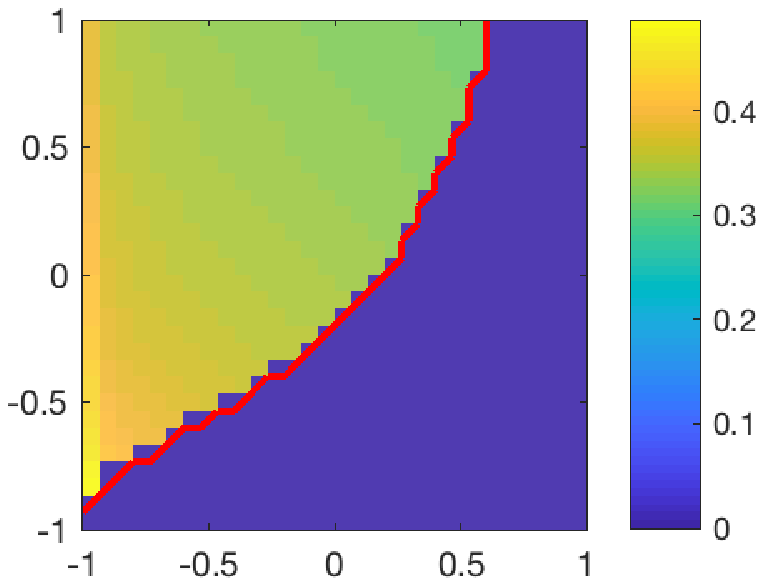}}
\subfigure[$\xi_3=1$.]
{\includegraphics[width=0.48\textwidth]{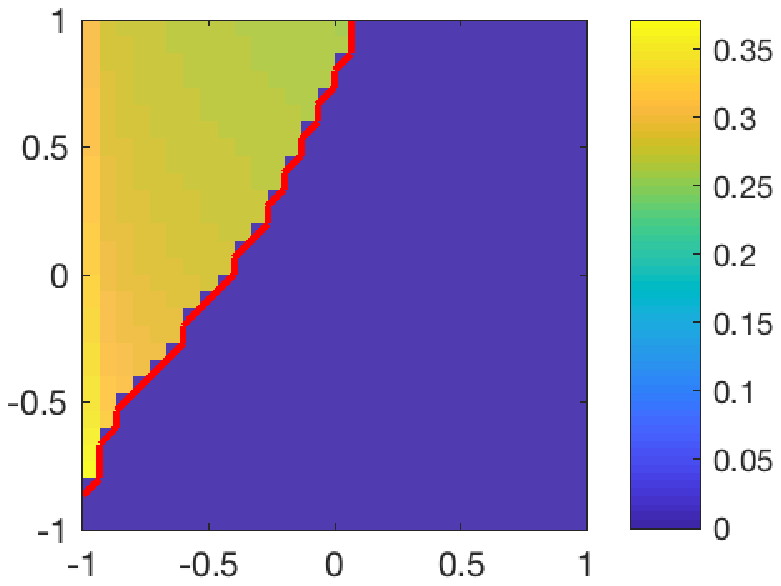}}
\caption{Exact solution $u(\xi_1,\xi_2, \xi_3, x, t)$ and computed zero level set $\{ \xi_1,\xi_2 | \phi = 0, \xi_3, x, t \}$ (solid red curve) at fixed $\xi_3=-1,0,1$, at the time $t=600$ (years) and location $x=600$, i.e., 100 m downstream from the injection point.}
	\label{fig:level_set_cuts}
\end{figure}

The level set problem is solved on a $31\times31\times 31$ grid in stochastic space until locking of the zero level set. The exact solution and the computed zero level set contour are depicted in Figure~\ref{fig:level_set_cuts} as a function of $\xi_1, \xi_2$ with $\xi_3$ fixed.

The zero level set agrees relatively well with the exact solution, but there is an error which becomes more severe close to the boundaries ($\xi_i=\pm 1$, $i=1,2,3$). It is also affected by the fact that the zeros of the speed function $F$ do not exactly coincide with the discontinuity of the solution $u$. Despite the alleged versatility of the level set method to accurately capture complex shapes, it was in this case difficult to numerically compute the zero level set.

Next, we use the level set solution to compute a spectral expansion. Both local polynomial basis functions on a simplex tessellation, and global frames will be used. In three (or more) dimensions, the tessellation becomes more challenging and it is nontrivial to make sure that each simplex element contains a sufficient number of solution evaluations. On the other hand, the global approach with only two solution regions determined by the sign of the level set function does not suffer from this issue.

\begin{figure}[H]
\centering
{\includegraphics[width=0.83\textwidth]{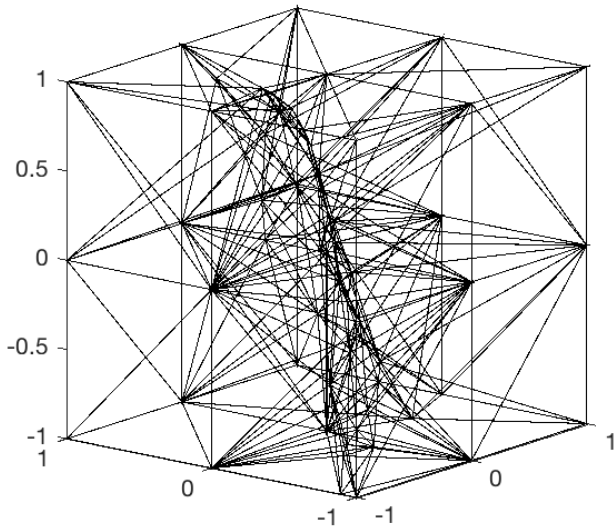}}
\caption{Simplex tessellation of the stochastic domain with 243 elements based on the computed zero level set.}
	\label{fig:3d_simplex_tess}
\end{figure}

For SOP, the stochastic domain is tessellated into simplices based on the computed zero level set, and a local orthonormal basis is computed for each simplex, as depicted in Figure~\ref{fig:3d_simplex_tess}. There is a tradeoff between a sufficiently fine simplex tessellation that conforms well with the zero level set, and one that is sufficiently coarse to ensure that a large enough number of solution evaluations are contained within each simplex for accurate computation of the local SOP coefficients. This issue is reflected in Table~\ref{tab:conv_diff_lev_3d}, where the relative error in the SOP solutions are shown for different simplex meshes, but with the same total number of solution evaluations. For some cases with a large number of simplex elements, the numerical error decreases with a lower-order polynomial reconstruction compared to a higher-order polynomial reconstruction. Of course, if the number of conservation law evaluations per simplex element is large enough, the error would decrease with increasing maximum polynomial order of the basis functions.
The number of basis functions per simplex element should remain below the number of solution evaluations in the same element. In a small number of elements this requirement is violated. Reconstruction of the stochastic solution is then performed through a reduced singular value decomposition.

\begin{table}[h]
\begin{tabular}{c|cc|cc|cc|}
\cline{2-7}
& \multicolumn{2}{ c| }{$N=1$, $P=4$} & \multicolumn{2}{ c| }{$N=2$, $P=10$} & \multicolumn{2}{ c| }{$N=3$, $P=20$}  \\ 
\cline{1-7}
\multicolumn{1}{|c|}{ $|S_{\be}|$ } & LAD & OLS & LAD & OLS & LAD & OLS  \\
\cline{1-7}
\multicolumn{1}{|r|}{222} & 1.99e-1 & 2.33e-1 & 1.61e-1 & 1.85e-1 & 1.38e-1 & 1.54e-1  \\
\multicolumn{1}{|r|}{675} & 1.40e-1 & 1.91e-1 & 1.26e-1 & 1.48e-1 & 1.17e-1 & 1.35e-1  \\
\multicolumn{1}{|r|}{1238} & 8.78e-2 & 1.13e-1 & 8.10e-2 & 9.63e-2 & 9.54e-2 & 1.07e-1  \\
\multicolumn{1}{|r|}{2047} & 6.04e-2 & 7.78e-2 & 7.60e-2 & 8.67e-2 & 1.27e-1 & 1.34e-1  \\
\cline{1-7}
\end{tabular}
\caption{Numerical convergence of $\epsilon^{\scriptscriptstyle{\textup{SOP}}}_{\ell_{1}}$ for the CO$_2$ transport problem for different simplex tessellations of the stochastic domain, and different orders of piecewise polynomial reconstruction. Ordinary Least Squares (OLS) and Least Absolute Deviations (LAD) are used to locally estimate the SOP coefficients.}
\label{tab:conv_diff_lev_3d}
\end{table}

\begin{figure}[H]
\centering
\subfigure[Reconstruction using LAD.]
{\includegraphics[width=0.48\textwidth]{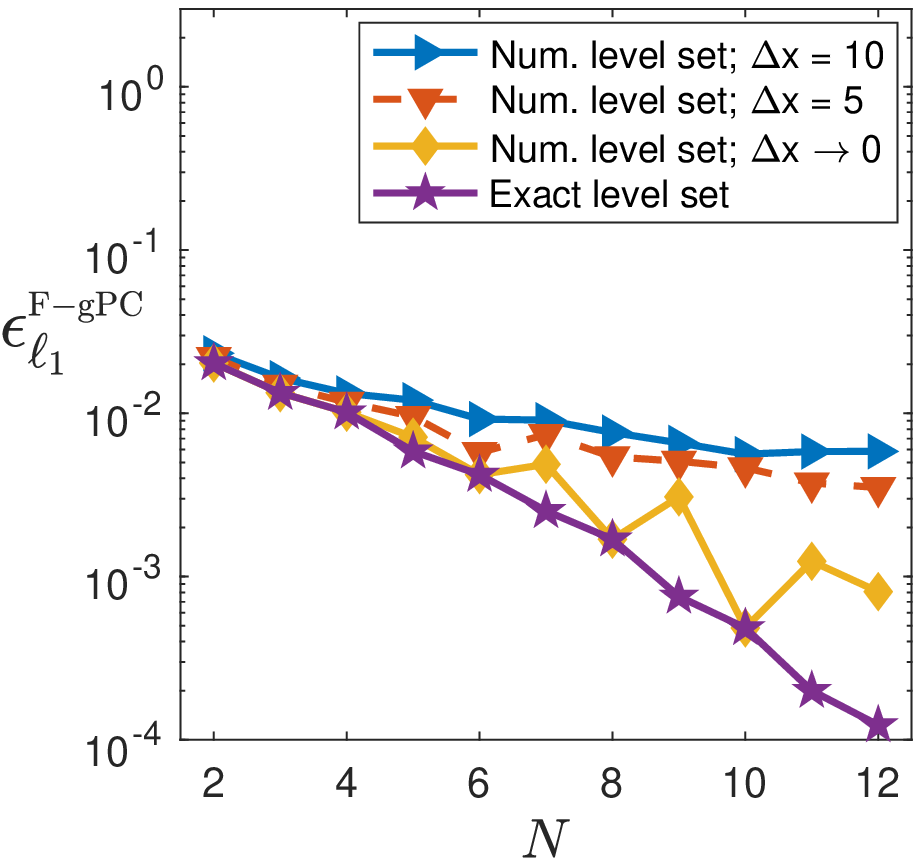}}
\hspace{0.1cm}
\subfigure[Reconstruction using OLS.]
{\includegraphics[width=0.48\textwidth]{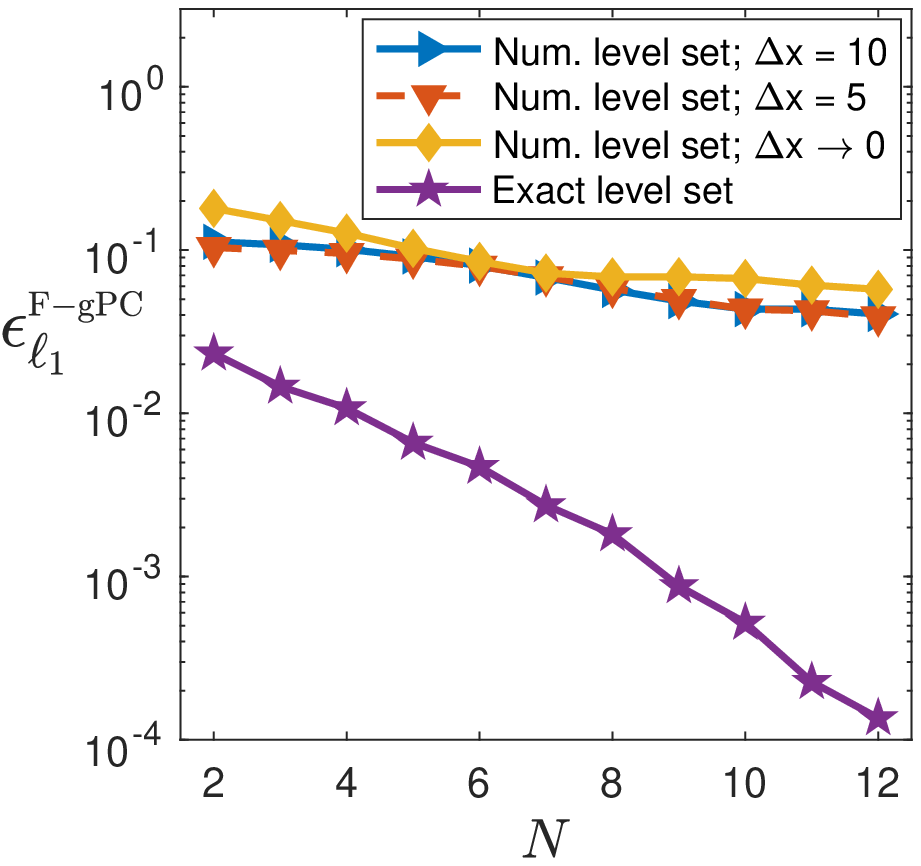}}

\caption{Numerical convergence of $\epsilon^{\scriptscriptstyle{\textup{F-gPC}}}_{\ell_{1}}$ for the CO$_2$ transport problem for different orders of piecewise polynomial frames, using OLS and LAD. Global total order Legendre polynomials are used for the frames, on a stochastic grid with 21 points per dimension.}
	\label{fig:conv_global_pol_CO2_l1_m21}
\end{figure}
\begin{figure}[H]
\centering
\subfigure[Reconstruction using LAD.]
{\includegraphics[width=0.48\textwidth]{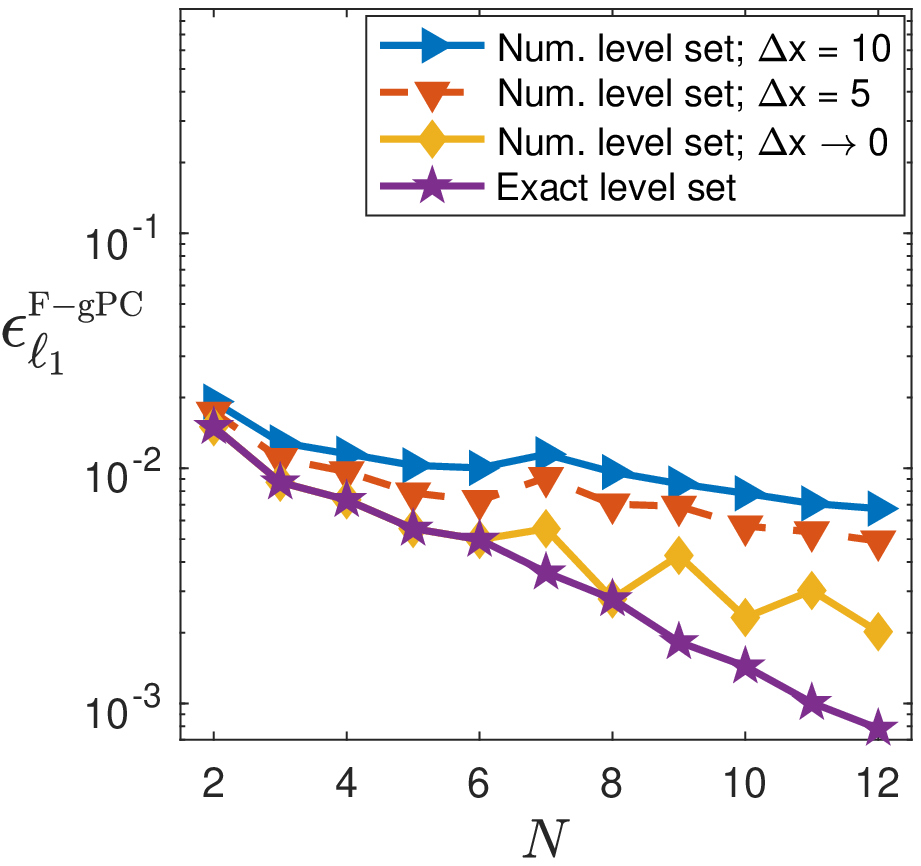}}
\hspace{0.1cm}
\subfigure[Reconstruction using OLS.]
{\includegraphics[width=0.48\textwidth]{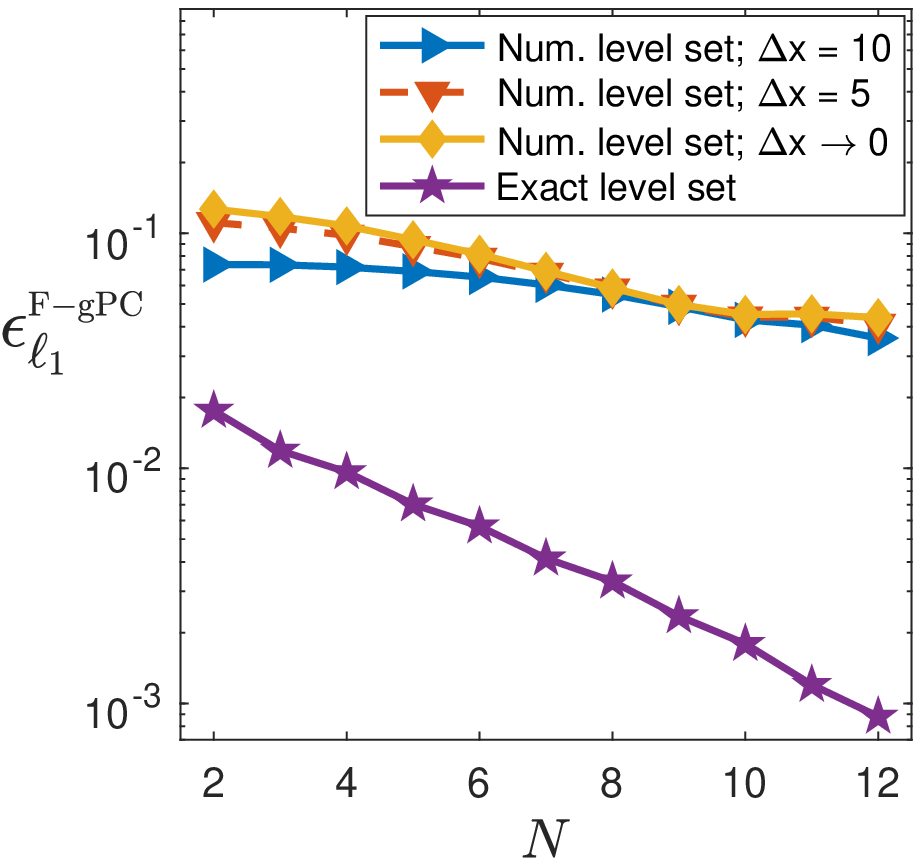}}

\caption{Numerical convergence of $\epsilon^{\scriptscriptstyle{\textup{F-gPC}}}_{\ell_{1}}$ for the CO$_2$ transport problem for different orders of piecewise polynomial frames, using OLS and LAD. Global total order Legendre polynomials are used for the frames, on a stochastic grid with 41 points per dimension.}
	\label{fig:conv_global_pol_CO2_l1_m41}
\end{figure}

Next, we use the computed level set solution to define frames by restricting orthogonal Legendre polynomials to the ranges of $\bxi$ where the level set function is negative and positive, respectively. Figures~\ref{fig:conv_global_pol_CO2_l1_m21} and~\ref{fig:conv_global_pol_CO2_l1_m41} show the relative $\ell_1$ errors in the solution of the CO$_2$ transport problem with F-gPC, varying the number of frame functions per solution region ($P$) as a function of maximum order of polynomial degree $N$. The level set problem is discretized with 21 and 41 points per dimension, respectively. In addition to the semi-analytical solutions, we present numerical results where the CO$_2$ transport problem has been solved with a second-order flux-limited Kurganov-Tadmor scheme~\cite{Kurganov_Tadmor_90}.

\begin{figure}[H]
\centering
\subfigure[Reconstruction using LAD.]
{\includegraphics[width=0.48\textwidth]{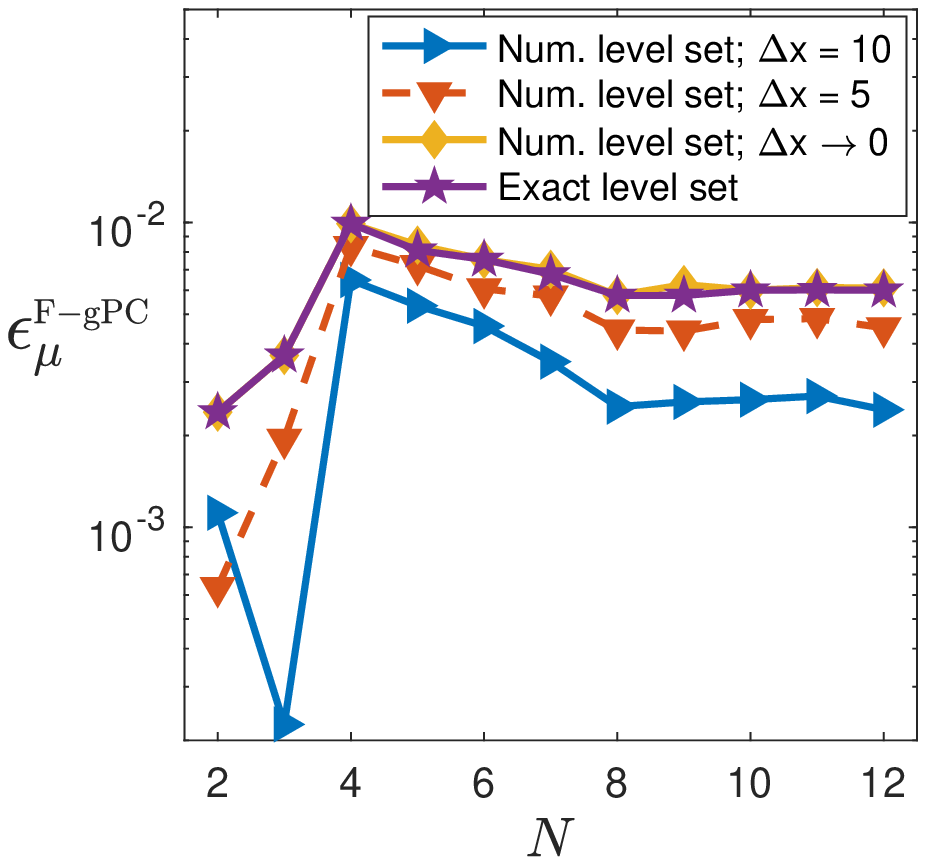}}
\hspace{0.1cm}
\subfigure[Reconstruction using OLS.]
{\includegraphics[width=0.48\textwidth]{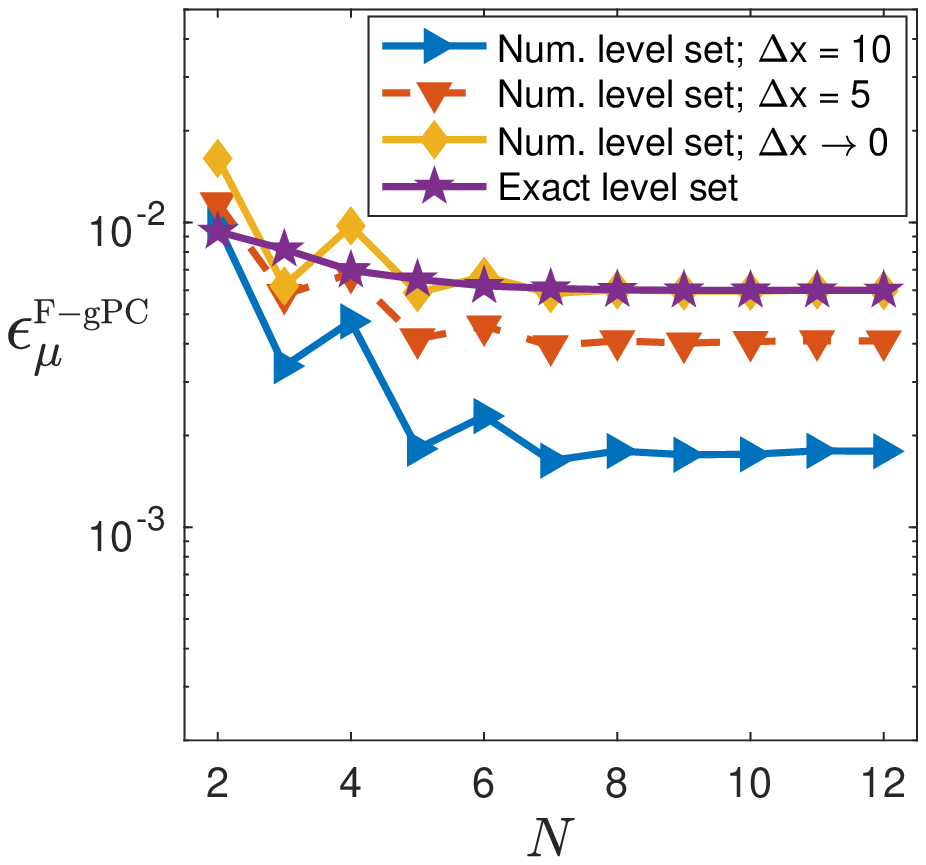}}
\subfigure[Reconstruction using LAD.]
{\includegraphics[width=0.48\textwidth]{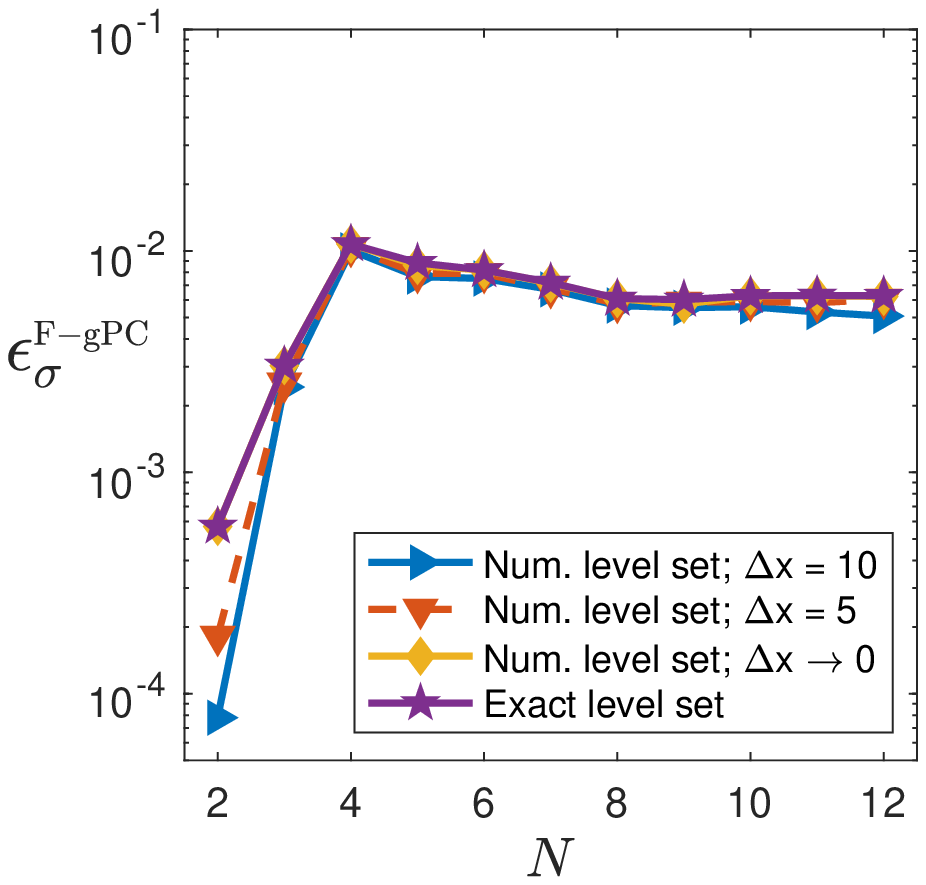}}
\hspace{0.1cm}
\subfigure[Reconstruction using OLS.]
{\includegraphics[width=0.48\textwidth]{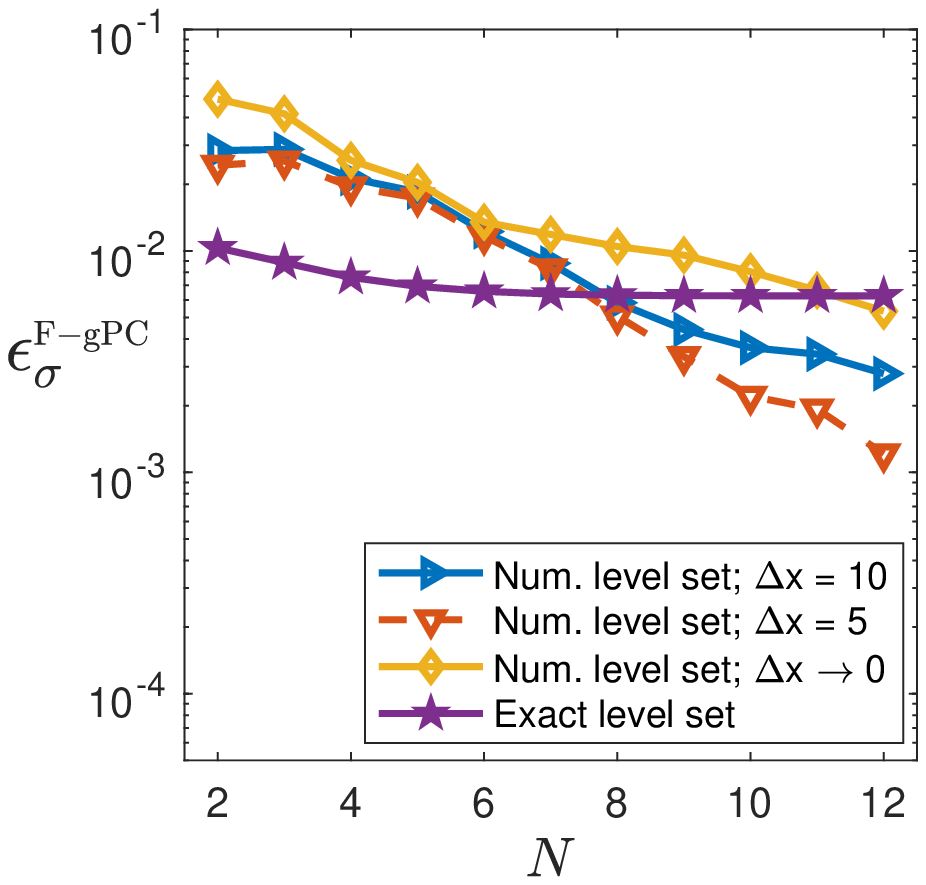}}

\caption{Numerical convergence of relative error in mean and standard deviation for the CO$_2$ transport problem for different orders of piecewise polynomial frames, using OLS and LAD. Global total order Legendre polynomials are used for the frames, on a stochastic grid with 21 points per dimension.}
	\label{fig:conv_global_pol_CO2_l1_m21_mean_std}
\end{figure}

\begin{figure}[H]
\centering
\subfigure[Reconstruction using LAD.]
{\includegraphics[width=0.48\textwidth]{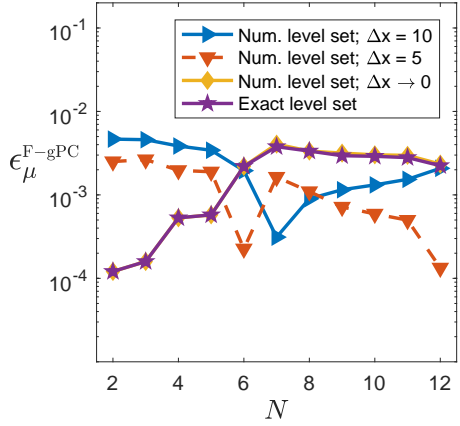}}
\hspace{0.1cm}
\subfigure[Reconstruction using OLS.]
{\includegraphics[width=0.48\textwidth]{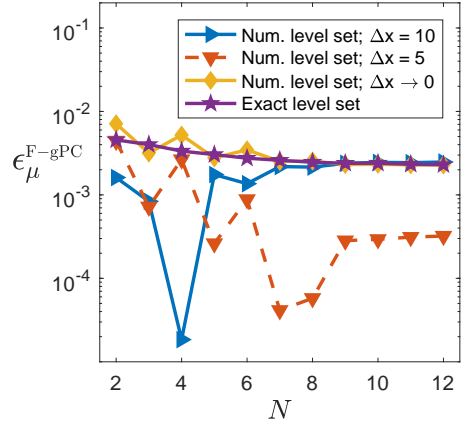}}
\subfigure[Reconstruction using LAD.]
{\includegraphics[width=0.48\textwidth]{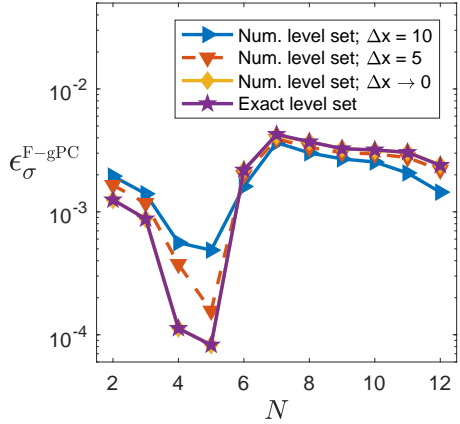}}
\hspace{0.1cm}
\subfigure[Reconstruction using OLS.]
{\includegraphics[width=0.48\textwidth]{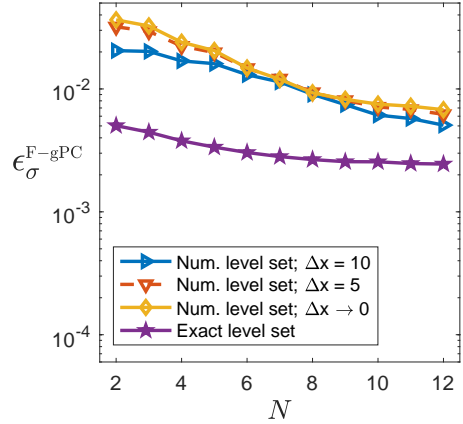}}

\caption{Numerical convergence of relative error in mean and standard deviation for the CO$_2$ transport problem for different orders of piecewise polynomial frames, using OLS and LAD. Global total order Legendre polynomials are used for the frames, on a stochastic grid with 41 points per dimension.}
	\label{fig:conv_global_pol_CO2_l1_m41_mean_std}
\end{figure}

Qualitatively, the behavior of the error is similar to the error using F-gPC for the 2D test case. For LAD, the total error is dominated by the discretization error of the CO$_2$ transport problem. The OLS error is around 1-3 orders of magnitude larger than the reference error. Again, this error is dominated by misclassified conservation law evaluations that are partially resolved by the polynomial reconstruction. This explains why the errors shown in Figures~\ref{fig:conv_global_pol_CO2_l1_m21} and~\ref{fig:conv_global_pol_CO2_l1_m41} do not exhibit monotone convergence with increasing polynomial order.

For completeness, the relative errors in means and standard deviations are included in Figures~\ref{fig:conv_global_pol_CO2_l1_m21_mean_std} and~\ref{fig:conv_global_pol_CO2_l1_m41_mean_std}, respectively. Similar to the previous test case, the error appears to be dominated by the numerical approximation error in these statistics. With improved approximation of mean and standard deviation using frames, we expect a similar trend as for the relative $\ell_1$ error.

For the numerical level set solution, the relative error is orders of magnitude smaller than the same error using the exact level set location, when using OLS. This demonstrates that accurate computation of the level set function is important. The smaller relative error on the coarse grid compared to the fine grid for the highest order expansions may seem surprising. The explanation is that the error is dominated by polynomial truncation and it is itself decribed by a highly oscillatory polynomial function. When poorly resolved, it may result in a smaller error on the coarser grids. The error using $\ell_1$ regression (LAD) is significantly smaller than the error using $\ell_2$ regression (OLS) for the numerical level set function, reflecting that it is less sensitive to 'outliers', i.e., conservation law evaluations from the other side of the zero level set.

The number of conservation law evaluations for ME-gPC is up to 200 times higher than the number of conservation law evaluations using the $21\times 21 \times 21$ grid discretization of stochastic space used to obtain similar errors shown in Figure~\ref{fig:conv_global_pol_CO2_l1_m21}. If the accuracy of the level set solver is improved, the gain would be even higher as shown for the exact zero level set results, also in Figure~\ref{fig:conv_global_pol_CO2_l1_m21}. If the level set is known to sufficient accuracy, the error would be 100 times smaller than for the adaptive ME-gPC method, in addition to the already significantly reduced computational cost. 

The reduced number of conservation law solves for the level set formulation has to be compared with the extra cost of solving the level set problem. For many complex problems, the numerical cost of the solution of the level set problem is small in comparison to that of a large number of calls to the solver of the conservation law, and a reliable proxy for the total numerical cost is then given by the total number of conservation law evaluations.


\section{Conclusions}
\label{sec:conc}
We have introduced a level set method to track discontinuities in the solutions of conservation laws in stochastic space by solving a Hamilton-Jacobi equation with a speed function that vanishes at discontinuities. The method is an adaptive surrogate method in the sense that the level set problem is solved on a sequence of successively finer grids in stochastic space, and high-fidelity conservation law solutions are replaced by interpolated estimates in regions of smoothness.

The level set solution can be used in various ways to reconstruct a proxy of the solution of interest to be used in fast postprocessing to obtain QI. A simplex tessellation of the stochastic domain leads to localized support of the stochastic basis functions and a set of small independent regression problems for the local simplex basis coefficients. While this in principle leads to more robustness with respect to misalignment of the computed level set with the exact discontinuities, the tessellation itself leads to numerical errors. For the problem setups investigated in this paper, frames perform better than simplex tessellations. For problems where the discontinuity has an end-point within the stochastic domain, e.g., the solution to the Kraichnan-Orszag problem, a simplex tessellation aligned with the discontinuity probably yields better results than a global gPC approximation. For such problems, frames cannot be used in a straightforward manner.

Significantly reduced computational cost, as measured by the total number of calls to the conservation law solver, has been demonstrated for the level set method with frame based solution reconstruction compared to adaptive ME-gPC. If the zero of the level set function is known exactly or to sufficient accuracy, the decay of the relative error in the stochastic solution is fast with respect to increasing polynomial order. We have observed up to two orders of magnitude smaller error compared to adaptive ME-gPC, in addition to a reduced cost from a smaller number of conservation law solutions. This is in contrast to the case of numerically computed level set which introduces an additional error, although the numerical cost comparison is still favorable compared to adaptive ME-gPC. This demonstrates that accurate solution and perhaps even an improved formulation of the level set problem is of great interest. The use of adaptive and unstructured grids for improving the discretization of the level set equation is an important direction for future work.


\section*{Acknowledgements}
Per Pettersson was supported by the Research Council of Norway through the project 244035/E20 CONQUER, under the CLIMIT program.

This material is based upon work of Alireza Doostan supported by the U.S. Department of Energy Office of Science, Office of Advanced Scientific Computing Research, under Award Number DE-SC0006402 and NSF grant CMMI-1454601.

\end{document}